\documentclass[
11pt, 
a4paper, 
oneside, 
]{article}
\usepackage[dvipsnames]{xcolor}
\usepackage{hyperref}
\usepackage[OT2,T1]{fontenc} 
\usepackage[light]{kpfonts} 
\usepackage[ left=2.5cm, right=2.5cm,top=4cm,bottom=4cm,]{geometry}
\usepackage[utf8]{inputenc} 
\usepackage{enumitem}
\usepackage{graphicx} 
\graphicspath{{Figures/}} 
\usepackage{amsmath,amssymb,amsthm,bbm} 
\usepackage{mathtools} 
\usepackage{varioref} 
\usepackage{tikz-cd} 
\usepackage[title]{appendix}
\usepackage{marginnote}
\usepackage{xspace}
\usepackage{bbm}
\usetikzlibrary{trees}

\newcommand{\Z}{\ensuremath{\mathbb{Z}}\xspace}
\newcommand{\Q}{\ensuremath{\mathbb{Q}}\xspace}
\newcommand{\R}{\ensuremath{\mathbb{R}}\xspace}
\newcommand{\C}{\ensuremath{\mathbb{C}}\xspace}

\newcommand{\F}{\ensuremath{\mathbb{F}}\xspace}

\newcommand{\comment}[1]{}

\DeclareMathOperator{\Gal}{Gal}

\DeclareMathOperator{\Hom}{Hom}

\newcommand{\hol}{\mathrm{hol}}

\newcommand{\GL}{\ensuremath{\mathrm{GL}}\xspace}
\newcommand{\SL}{\ensuremath{\mathrm{SL}}\xspace}

\newcommand{\af}{\mathfrak{a}\xspace}

\newcommand{\mb}{\mathbb}
\newcommand{\mc}{\mathcal}

\newcommand{\oo}{\mathcal{O}}

\newcommand{\PGL}{\mathrm{PGL}}

\newcommand{\cG}{\mathcal{G}}

\usepackage{color}
\usepackage{todonotes}
\setuptodonotes{color=yellow!50}
\usepackage{comment}
\usepackage[normalem]{ulem}

\theoremstyle{definition} 
\newtheorem{definition}{Definition}
\theoremstyle{plain} 
\newtheorem{theorem}{Theorem}[subsection]

\newtheorem{corollary}{Corollary}
\newtheorem{conjecture}{Conjecture}

\theoremstyle{remark} 
\newtheorem{remark}{Remark}
\newtheorem{example}{Example}

\input xy
\xyoption{all}

\DeclareMathOperator{\disc}{disc}
\DeclareMathOperator{\Nm}{Nm}

\newcommand{\cO}{\mathcal {O}}

\newcommand{\cH}{\mathcal{H}}

\newcommand{\cW}{\mathcal{W}}

\newcommand{\ST}{\operatorname{ST}^{\times}}
\newcommand{\ord}{\operatorname{ord}}

\newcommand{\mbP}{\mathbb{P}}
\DeclareMathOperator{\red}{red}

\DeclareMathOperator{\order}{ord}
\DeclareSymbolFont{cyrletters}{OT2}{wncyr}{m}{n}
\DeclareMathSymbol{\Sha}{\mathalpha}{cyrletters}{"58}
\newcommand{\Hc}{\ensuremath{\mathcal{H}}\xspace}

\DeclareMathOperator{\RM}{\mathrm{RM}}
\DeclareMathOperator{\triv}{triv}
\DeclareMathOperator{\res}{res}

\hypersetup{linktocpage,bookmarksopen = true, bookmarksopenlevel = 1, colorlinks=true,
linkcolor = Bittersweet, citecolor = RoyalBlue , urlcolor = gray, pdfstartview = FitR}

\begin{document}
	
	\title{Real quadratic singular moduli and \\$p$-adic families of modular forms}
	\author{Paulina Fust, Judith Ludwig, Alice Pozzi,  Mafalda Santos and Hanneke Wiersema}
	\date{}

	\maketitle
	
\begin{abstract}
The classical theory of elliptic curves with complex multiplication is a fundamental tool for studying the arithmetic of abelian extensions of imaginary quadratic fields. 
While no direct analogue is available for real quadratic fields, a (conjectural) theory of ``real multiplication'' was recently proposed by Darmon and Vonk, relying on $p$-adic methods, and in particular on the new notion of rigid meromorphic cocycles.  
A rigid meromorphic cocycle is a class in the first cohomology of the group $\mathrm {SL}_2(\mathbb Z[1/p])$ acting on the non-zero rigid meromorphic functions on the Drinfeld $p$-adic upper half plane by M\"obius transformation. The values of rigid meromorphic cocycles at real quadratic points can be thought of as analogues of singular moduli for real quadratic fields.

In this survey article, we will discuss aspects of the theory of complex multiplication and compare them with conjectural analogues for real quadratic fields, with an emphasis on the role played by families of modular forms in both settings.
\end{abstract}

	\section{Introduction}
The goal of this article is to present recent developments in the theory of \emph{singular moduli for real quadratic fields} introduced in \cite{DV2021} in a parallel perspective to the classical theory of complex multiplication. 
The theory of elliptic curves with complex multiplication has yielded some spectacular arithmetic applications. Kronecker tackled the problem  of constructing abelian extensions of imaginary quadratic fields, known as the \emph{Kronecker Jugendtraum}, via singular moduli, i.e., values of modular functions at imaginary quadratic points of the complex upper half plane $\cH_{\infty}$. While class field theory provides a description of the Galois groups classifying abelian extensions of arbitrary number fields, finding an explicit recipe for their generators is more elusive beyond the CM setting.

	Another problem for which in current approaches CM theory is an essential tool is the Birch and Swinnerton-Dyer Conjecture. 
	Given an elliptic curve $E$ over a number field $K$, the Birch and Swinnerton-Dyer Conjecture (BSD) predicts an equality between the rank of  its group of $K$-rational points and the order of vanishing of its $L$-function $L(E/K,s)$ at $s=1$. The only known approach to systematically producing global points is via the theory of Heegner points. The works of Gross--Zagier \cite{GZ1985} and  Kolyvagin \cite{Kolyvagin}
	 settle many instances of BSD in rank 1 by exploiting properties of this supply of global points.
		
	The construction of singular moduli and Heegner points share some common features. They both arise from an appropriate evaluation process at CM points on the complex upper half plane. The corresponding values are defined over abelian extensions of imaginary quadratic fields $K$. We will refer to this type of construction as \emph{Heegner constructions}. It is natural to ask if a similar approach can be implemented to construct collections of objects lying in abelian extensions of other number fields.  
	
	Let $K$ be a real quadratic field. From a naive perspective, the idea of evaluating modular functions at real quadratic points (or RM, by analogy with CM points) of $\cH_{\infty}$ already fails because RM points lie on the real line. This suggests that an appropriate analogue of CM theory should be sought by replacing the role played by the archimedean place with a finite place that does not split in $K$. This is the point of view proposed by Darmon and Vonk in their theory of singular moduli for real quadratic fields.

	For a finite prime $p$, the Drinfeld $p$-adic upper half plane is the rigid analytic space $\cH_p$ whose $\C_p$-valued points are $\mbP^1(\C_p) \setminus \mbP^1(\Q_p)$. It contains RM points over real quadratic fields $K$ in which $p$ is ramified or inert. In analogy with CM theory, it would be natural to attempt to define singular moduli as values of suitable modular functions at RM points in $\cH_p$, for some action of a modular group $\Gamma \subset \mathrm{GL}_2(\Q_p). $ Typically, one would expect this group to arise as the elements of norm one in a $\Z\left [1/p \right ]$-order of a quaternion algebra $B/\Q$ which splits at $p$. On the one hand, if the quaternion algebra is definite, the quotient $\Gamma \backslash \cH_p$ admits the structure of a rigid analytic space. However, in that setting, $B$ does not contain any real quadratic subalgebra. Therefore one would not expect such quotients to appear in the framework of RM theory. On the other hand, if $B$ is indefinite, it contains real quadratic subalgebras but the action of $\Gamma$ on $\cH_p$ is not discrete, so that there are no suitable modular functions. 	Darmon and Vonk fix the indefinite quaternion algebra $B=\mathrm{M}_2(\Q_p)$, so that $\Gamma=\SL_2(\Z[1/p])$ is the Ihara group, and they define rigid meromorphic cocycles as classes in $H^1(\Gamma, \mathcal {M}^{\times}),$ where $\mathcal M$ is the ring of rigid meromorphic functions on $\cH_p$. They  show that  these classes (or more precisely, the quasi-parabolic classes, see \S \ref{SS:  RMC}) 
	can be meaningfully evaluated at RM points.  They then conjecture that their RM values lie in composita of abelian extensions of real quadratic fields, making them suitable candidates for an analogue of singular moduli in the real quadratic setting.  The theory of real quadratic singular moduli fits into a program launched by Darmon in \cite{Dar01} attempting to define conjectural analogues of Heegner constructions for real quadratic fields via $p$-adic methods. Most notably, the construction of a putative analogue of Heegner points, known as Stark--Heegner points, has since seen many generalisations.  While the theoretical evidence on the rationality of these points is scarce, the computational results are extensive.	
	
	In the last 50 years, new aspects of the theory of complex multiplication have emerged, in which automorphic methods are used substantially. An important theme is the study of modular generating series constructed out of CM cycles on Shimura curves, pioneered by the work of Gross and Zagier on the factorisation of differences of singular moduli \cite{GZ1985}. The crucial idea of relating generating series to derivatives of real analytic families of Eisenstein series culminated in the proof of many cases  of the Birch and Swinnerton--Dyer Conjecture for a modular elliptic curve in analytic rank 1 in \cite{GZ1986} (see Thm. \ref{T: GZ Inventiones}). Although the algebraicity properties of the analogues of the Heegner constructions in the real quadratic setting are conjectural, one can package RM cycles on $\mathcal H_{p}$ into modular generating series. It is natural to speculate that these generating series would be related to  derivatives of $p$-adic families of modular forms. This perspective has been explored by Darmon, Pozzi and Vonk in \cite{DPV1} and \cite{DPV2}.
	
	In the archimedean setting, only real analytic Eisenstein series  admit variations in families. By contrast, the work of Hida, and its generalisation via the theory of eigenvarieties introduced by Coleman--Mazur, show that $p$-adic families of automorphic forms abound. A useful feature of the $p$-adic theory is that cusp forms vary in $p$-adic families. In addition, their corresponding Galois representations can be interpolated.  This aspect has been leveraged into breakthroughs in Iwasawa theory and the theory of $p$-adic $L$-functions, for example the proof of the Iwasawa Main Conjecture over totally real fields by Mazur--Wiles, or the proof of the  Mazur--Tate--Teitelbaum Conjecture by 	Greenberg--Stevens.
	
	 In the context of modular generating series for RM values, the additional flexibility of the $p$-adic setting has been used by Darmon--Pozzi--Vonk to prove some instances of the expected algebraicity results for the analogues of elliptic units for real quadratic fields (see Theorem \ref{T: DPV2}). 
Their approach is somewhat indirect. On the one hand, derivatives of certain $p$-adic families of modular forms are used to produce modular generating series for RM values. On the other hand, the derivatives of cuspidal families are related to global cohomology classes informing the arithmetic of abelian extensions of real quadratic fields. 
	This strategy towards the algebraicity of Heegner constructions in the RM setting suggests that the richer theory of $p$-adic automorphic forms can sometimes make up for the lack of an \emph{a priori} global description of the corresponding RM objects.  \\

\par 
\textbf{Overview of the paper.} This article is structured as follows. Section \ref{S: Heegner CM} contains a review of some aspects of the classical theory of complex multiplication, with an emphasis on singular moduli, elliptic units and Heegner points. Section \ref{S: Heegner RM} focuses on the construction of their conjectural analogues for real quadratic fields. For the reader's convenience, we include some background on the Drinfeld $p$-adic upper half plane as a rigid space and  introduce the theory of rigid meromorphic cocycles. 
	Sections \ref{S: analytic families} and \ref{S: p-adic families} are devoted to the interplay between derivatives of families of modular forms and singular moduli, in the archimedean and non-archimedean setting respectively. In particular, Section \ref{S: analytic families} treats Gross and Zagier's work on the factorisation of singular moduli and the BSD Conjecture in analytic rank 1 \cite{GZ1985} and \cite{GZ1986}, while Section \ref{S: p-adic families} focuses on the articles \cite{DPV1} and \cite{DPV2}.
	Sections \ref{S: Heegner CM} and \ref{S: Heegner RM} are meant to be mostly self-contained, while Sections \ref{S: analytic families} and \ref{S: p-adic families} aim to be an overview of various arithmetic applications to Heegner constructions in a unified framework. \\
	
		\par
\textbf{Acknowledgements.} We would like to thank the organisers of the WIN-Europe 4 conference for inviting us to contribute to these proceedings. We would like to thank Henri Darmon and Jan Vonk for enlightening explanations of their work. We would also like to thank Sandra Rozensztajn, Ana Caraiani and Catherine Hsu for helpful conversations during the preparation of this manuscript. We would like to thank O\u guz Gezmi\c s and the anonymous referees for many helpful comments and suggestions.\\
P.F. was funded by the DFG Graduiertenkolleg 2553.
J.L.\ acknowledges support from the Deutsche Forschungsgemeinschaft (DFG, German Research Foundation) through TRR 326 \textit{Geometry and Arithmetic of Uniformized Structures},  project number 444845124. A.P. was funded by the Leverhulme Trust (via the Leverhulme Prize of Ana Caraiani).  M.S. was funded by the Engineering and Physical Sciences Research Council [EP/L015234/1], the EPSRC Centre for Doctoral Training in Geometry
and Number Theory (The London School of Geometry and Number Theory), Imperial College
London, and by the Max-Planck-Institut für Mathematik. 
H.W. acknowledges support from the Herchel Smith Postdoctoral Fellowship Fund, and the Engineering and Physical Sciences Research Council (EPSRC) grant EP/W001683/1.

\section{Heegner Constructions for imaginary quadratic fields}
\label{S: Heegner CM}
In this section, we will give a brief overview of results in the classical theory of elliptic curves with complex multiplication, and its applications to the construction of arithmetic objects over abelian extensions of imaginary quadratic fields. 
At the heart of these results lies the fact that modular curves, quotients of the complex upper half plane by congruence subgroups, admit a moduli interpretation which endows them with a rational structure. The moduli interpretation of CM points can be used to determine their images under field automorphisms, which allows one to deduce that they define global points on modular curves.\\
We will begin by recalling some facts about the moduli interpretation of CM points on modular curves and the Main Theorem of Complex Multiplication. We will then discuss in more detail the Heegner constructions of singular moduli, elliptic units and Heegner points.

\subsection{Modular curves and CM points}
\subsubsection{Elliptic curves with complex multiplication}
Let $E$ be a complex elliptic curve, which we view as a Riemann surface. It admits a uniformisation $E(\C) \simeq \C/\Lambda$, where $\Lambda \subset \C$ is its period lattice, well-defined up to homothety; by rescaling, the latter can be chosen of the form $\Lambda =\tau \Z+ \Z$, for $\tau \in \C$ with positive imaginary part. Endomorphisms of $E$ are given by endomorphisms of its universal covering $\C$ that preserve the lattice $\Lambda$. Denote by $\mathcal R_{E}$ the set of endomorphisms of $E$. For a generic $\tau$ in the upper half plane, all endomorphisms of $E$ are given by multiplication by an integer. Otherwise, there exists $\alpha \in \C \setminus \Z$, and integers $a, b,c,d $ satisfying
\[
c\tau+d=\alpha, \qquad  a\tau+b=\alpha \tau. 
\]  
This forces $K=\Q(\tau)$ to be an imaginary quadratic field and $\alpha$ to be an algebraic integer in $K$. As a result $\mathcal R_E$ is isomorphic to an order $\mathcal O$ in $K$, for which the lattice  $\Lambda$ is a proper fractional ideal. In this case we say that $E$ has complex multiplication by $\mathcal O$. \\

\subsubsection{The Main Theorem of Complex Multiplication}
\par Let $\mathcal H_{\infty}$ be the Poincar\'e complex upper half plane, on which the group  $\GL^+_2(\R)$ of invertible real matrices with positive determinant acts  by Möbius transformation. \\
For a positive integer $N$, consider the order  
\[
\mathrm{M}_0(N)=
\left \{ \left (\begin{matrix}
	a       & b \\
	c       & d
	\end{matrix} \right ) \in \mathrm{M}_2(\Z): N\mid c  \right \}
\]
in the ring of $2\times 2$-matrices with integer coefficients, and let $\Gamma_0(N)$ be the group of matrices of determinant 1 in $M_0(N)$. The assigment $
\tau \mapsto (\tau \Z+ \Z, \tau \Z+1/N \Z)$  identifies the quotient $\Gamma_0(N) \backslash \cH_{\infty}$ with the moduli space of pairs of lattices 
$(\Lambda , \Lambda')$ with 
\[
\Lambda \subset \Lambda'  \text{ such that } \Lambda'/\Lambda \simeq \Z/N\Z, 
\]
 defined up to homothety.
  
  The moduli problem parametrising isomorphism classes of pairs of elliptic curves $(E, E')$ with a cyclic degree $N$-isogeny $E \to E'$ admits a coarse moduli space $Y_0(N)$ over $\Q$ which we refer to as the (open) modular curve of the level $\Gamma_0(N)$.
We let $X_0(N)$ be its compactification, obtained by adding cusps.

 There is a bijection  
\[
\iota_N \colon \Gamma_0(N) \backslash \mathcal H_{\infty} \to Y_0(N)(\C), \qquad 
\tau \mapsto (\C/(\tau \Z+ \Z), \ \C/(\tau \Z+1/N \Z))
\]
which provides a uniformisation of $Y_0(N)(\C)$ by the complex upper half plane. \\

We say that $\tau \in \mathcal H_{\infty}$ is imaginary quadratic of discriminant $D$ if $[\tau,1]$ is a solution to a quadratic equation of the form $A X^2+BXY +CY^2=0$ for some coprime integers $A$, $B$, $C$ and $D=B^2-4AC<0$. 
For $\tau \in \cH_{\infty}$ and $N\in \mathbb{Z}_{\geq 1}$ let
\[\mc{O}_\tau^{(N)} := \left \{ \begin{pmatrix}
	a       & b \\
	c       & d
	\end{pmatrix} \in \mathrm{M}_0(N) :  a \tau+b=c\tau^2+d\tau \right \}.
	\]
The image of the inclusion  $\mc{O}_{\tau}^{(N)} \hookrightarrow \C$ given by $\left ( \begin{smallmatrix} a & b\\ c &d  \end{smallmatrix} \right ) \mapsto c\tau+d$ is either isomorphic to $\Z$ or to an order in the imaginary quadratic field  $K=\Q(\tau)$. In the latter case, we say that $\tau$ is a CM point on $Y_0(N)$ and denote by $[\tau]$ its $\Gamma_0(N)$-orbit, when the level $\Gamma_0(N)$ is understood.
For an order $\mathcal O$ in an imaginary quadratic field $K$ of discriminant prime to $N$, we denote by $\mathcal H_{\infty}^{\mathcal O}$ the set of points $\tau \in \cH_{\infty}$ with 
\[
\mathcal O_{\tau}^{(N)}=\mathcal O_{\tau}^{(1)} \cap \mathcal O_{N\tau}^{(1)} \simeq \mathcal O.
\]  

Such points are in bijection with the pairs of elliptic curves $(\C/\Lambda, \C/\Lambda')$ whose corresponding lattices $(\Lambda, \Lambda')$ are contained in $K$ (up to homothety) and for  which the set 
\begin{equation} \label{E: proper}
\{ \alpha  \in K : \alpha \Lambda \subset \Lambda,\  \alpha \Lambda'\subset \Lambda' \}
\end{equation}
is equal to $\mathcal O$.
The set $\mathcal H_{\infty}^{\mathcal O}$ is preserved by the action of $\Gamma_0(N)$ and the quotient $\Gamma_0(N) \backslash \mathcal H_{\infty}^{\mathcal O}$ is finite.
We let $\Gamma_0(N) \backslash \mathcal H_{\infty}^\mathrm{CM}$ be the union of $\Gamma_0(N)$-orbits of points with complex multiplication by $\oo$ for some imaginary quadratic order $\oo$. \\
\par 

Denote by $\mathrm{Cl}(\mathcal O)$ the class group of $\mathcal O$. Given a proper fractional ideal $\mathfrak a $ in $\mathcal O$ and a pair of lattices $(\Lambda, \Lambda')$, condition \eqref{E: proper} is stable under multiplication by $\mathfrak a$. As a result 
\[
\mathfrak a *(\C/\Lambda,\C/ \Lambda'):=(\C/(\mathfrak a*\Lambda), \C/(\mathfrak a*\Lambda'))
\] 
defines an action of the class group $\mathrm{Cl}(\mathcal O)$ on 
$\Gamma_0(N)\backslash \mathcal H_{\infty}^{\mathcal O}$. For $[\tau] \in \Gamma_0(N)\backslash \mathcal H_{\infty}^{\mathcal O}$, we write 
$\mathfrak a*\tau$ for $\mathfrak a *(\C/(\Z+\tau \Z),\C/ (1/N\Z+\tau \Z))$ by abuse of notation.\\

\par 
The Main Theorem of Complex Multiplication ensures that CM points on modular curves are defined over suitable ring class fields, and describes the Galois action on CM points. 
More precisely, let $\mathcal O$ be an order in an imaginary quadratic field $K$. The Artin reciprocity map yields an isomorphism 
	\begin{equation*} \label{Artin rec}
	\mathrm{rec} \colon \mathrm{Cl}(\mathcal O) \to \Gal (H/K)
	\end{equation*}
	where $H$ is the ring class field of $\cO$. The Main Theorem of Complex Multiplication can be stated as follows.
\begin{theorem} \label{T: main CM}
Let $\oo$ be an order of discriminant $D$ in an imaginary quadratic field $K$, and let $H$ be its corresponding ring class field. Let $N$ be a positive integer coprime to $D$. \\For every point $[\tau] \in \Gamma_0(N) \backslash \mathcal H_{\infty}^{\oo}$, we have:
\begin{enumerate}[label=(\roman*)]
\item  $\iota_N([\tau]) \in Y_0(N)(H)$,
\item  (Shimura reciprocity law): every $\mathfrak a \in \mathrm{Cl}(\mathcal O)$  satisfies $\iota_N (\mathfrak a * [\tau])=\mathrm{rec}(\mathfrak a)^{-1}\iota_N([\tau])$.
\end{enumerate}
\end{theorem}
For a reference see \cite[Cor.\ 11.37]{Cox1989}.
In \S \ref{SS: Heegner Objects CM} we will explain how Theorem \ref{T: main CM} can be exploited towards the construction of global units and global points on elliptic curves over ring class fields of imaginary quadratic fields.

\subsection{Heegner Constructions}\label{SS: Heegner Objects CM}

Theorem \ref{T: main CM} endows the modular curve $Y_0(N)$ with a systematic supply of global points.  
Given a group scheme $X/\Q$ together with a $\Q$-rational morphism $Y_0(N) \to X$, one obtains a collection of global points on $X$, satisfying similar properties to those of CM points on modular curves. We refer to constructions of this flavour as \emph{Heegner constructions.}

In the following examples, Heegner constructions are obtained from an a priori holomorphic map \[
 \Phi \colon \Gamma_0(N) \backslash \cH_{\infty} \to X(\C) \]
of Riemann surfaces that is shown to arise from a $\Q$-rational morphism of curves (an important example being the modular parametrisation of an elliptic curve). 
As a consequence of Theorem  \ref{T: main CM}, the following properties are satisfied:
\begin{enumerate}[label=(\roman*)]
	\item If $[\tau] \in \Gamma_0(N) \backslash \cH_{\infty}^{\mathcal O}$ then $\Phi([\tau]) \in X(H)$; 
	\item The action of the class group $\mathrm{Cl}(\mathcal O)$ is compatible with the Galois action on $X(H)$ under a Shimura reciprocity law.  
\end{enumerate}

Below we will give three examples of Heegner constructions. The following table
summarises how these examples fit into this perspective.
\begin{center}
\begin{tabular}{c|c|c}
 & $X$  & $\Phi$     \\[1pt] 
 \hline \hline   
 & \\[\dimexpr-\normalbaselineskip+3pt]
 Singular Moduli & affine line $\mathbb A^1/\Q$ & $j$-invariant      \\
 Elliptic Units & multiplicative group $\mathbb G_m/\Q$ &  modular unit    \\
 Heegner Points &  elliptic curve $E/\Q$  & modular parametrisation  
\end{tabular}
\end{center}
The setting of singular moduli is obtained by specialising Theorem~\ref{T: main CM} to the case of level $\Gamma_0(1)$. Elliptic units are constructed from the evaluation of a modular unit at CM points. Roughly speaking, Heegner points are obtained by projecting CM cycles on the Jacobian of the modular curve onto its $f$-isotypic component for a newform $f$ of weight 2. \begin{remark}
The constructions presented in \S \ref{S: Heegner RM} will mimic the ones above by replacing $\cH_{\infty}$ with the Drinfeld $p$-adic upper half plane $\cH_p$ and imaginary quadratic with real quadratic points. The Heegner constructions in the RM setting are obtained via techniques of rigid analytic geometry and are expected to satisfy similar properties to those of classical Heegner constructions. However, there are significant differences. The analogues of the map $\Phi$ described above will only be defined at real multiplication points. More crucially, the RM constructions are purely local, and no analogues of the modular curve are available in that setting.
\end{remark}

We will now proceed to describe the inputs of these constructions in some detail.

\subsubsection{Singular moduli}\label{SSS: singular moduli}
Given a lattice $\Lambda$ in $\C$, the theory of elliptic curves over the complex numbers yields an identification between the complex torus $\C/\Lambda$ and the points of an elliptic curve with equation 
\[
Y^2=4X^3-g_2(\Lambda)X-g_3(\Lambda)
\]
where the constants are defined as 
\[
g_2(\Lambda)=60 \sum_{\omega \in \Lambda-\{0\}} \frac{1}{\omega^4}, \quad g_3(\Lambda)=140 \sum_{\omega \in \Lambda-\{0\}} \frac{1}{\omega^6}. 
\]
For $i=2, 3$, and $\tau $ in the upper half plane, the assignment 
$g_i(\tau):=g_i(\Z+\tau \Z)$ defines a modular form of weight $2i$ and level $\SL_2(\Z)$. 
The complex functions 
\begin{align*}
\Delta:=g_2^3-27g_3^2,  && 
j:=1728 \frac{g_2^3}{\Delta}
\end{align*}
are the \emph{modular discriminant} and the \emph{modular $j$-function} respectively. The former is a nowhere vanishing cusp form of weight $12$. The latter is a holomorphic $\SL_2(\Z)$-invariant function on $\mathcal H_{\infty}$. The map  
\begin{equation} \label{E: j-function}
j\colon \SL_2(\Z) \backslash \mathcal H_{\infty}\to \mathbb A^1 (\C)
\end{equation}
extends to an isomorphism of Riemann surfaces between the compactifications of $\SL_2(\Z) \backslash \mathcal H_{\infty}$ and $\mathbb A^1(\C)$, respectively. 
The crucial property of the $j$-invariant is that it characterises elliptic curves up to isomorphism over an algebraically closed field, and captures the minimal field of definition of an elliptic curve. As a result, the coarse moduli space for elliptic curves  $Y_0(1)$ can be identified with the affine line via \eqref{E: j-function}. For $N\geq 1$, the image of the morphism 
\begin{equation*}
\Gamma_0(N) \backslash \mathcal H_{\infty}\to \mathbb A^2 (\C)\ , \qquad \tau \mapsto (j(\tau), j(N\tau))
\end{equation*}
is cut out by a single equation given by $\varphi_N \in \Q[X,Y]$, referred to as the $N$-th modular polynomial (\cite[Section 7.5]{DS}). 
The corresponding plane curve is birational to the modular curve $X_0(N)$. 

\begin{definition}
\emph{Singular moduli} are values of the $j$-function at points in $\SL_2(\Z) \backslash \mathcal H_{\infty}^{\mathrm{CM}}$. 
\end{definition}

For $N=1$, Theorem \ref{T: main CM} directly translates into properties of singular moduli. The general level case easily follows from the level 1 setting, together with properties of the modular polynomials.
An important application is the following corollary.
\begin{corollary}
Let $D<0$ be a fundamental discriminant, let $\tau \in \mathcal H_{\infty}$ be a point of discriminant $D$ and denote $K=\Q(\tau).$   Then the Hilbert class field of $K$ is $K(j(\tau)).$
\end{corollary}  
This result (and its analogues for non-fundamental discriminants) can be interpreted as providing explicit generators for certain abelian extensions of imaginary quadratic fields. More generally, Hilbert's Twelfth Problem consists in finding explicit formulae for generators of all abelian extensions of number fields. The motivating example is the Kronecker--Weber Theorem. It states that all abelian extensions of $\Q$ can be obtained by adjoining roots of unity, which can be viewed as values of the function $z \mapsto e^{2\pi i z}$ at rational arguments. However, for an imaginary quadratic field $K$, singular moduli do not suffice to generate the maximal abelian extension, as one would also need to consider values of coordinates of torsion points of elliptic curves with CM by orders in $K$ (see  \cite[Chp. II, Cor. 5.7]{Silverman}).

\subsubsection{Elliptic units}\label{SSS: elliptic units}
Singular moduli give generators of certain abelian extensions of imaginary quadratic fields. For arithmetic applications, it can be useful to construct generators with controlled integrality properties, by replacing the $j$-function with suitable units on modular curves of higher level.

A modular unit $u$ is a nowhere vanishing function on the open modular curve $Y_0(N)$ extending meromorphically to the cusps. Examples of modular units can be constructed as products of pullbacks of the modular discriminant via projection maps from $Y_0(N)$ to the modular curve of level 1. For a fixed level $N$, let $\Sigma_N$ be the set of divisors of $N$ and consider a degree zero formal linear combination $\alpha=\sum_{d \in \Sigma_N} m_d[d]$. Then the holomorphic function 
\[
u_{\alpha}(z):=\prod_d\Delta(dz)^{m_d}
\]
is $\Gamma_0(N)$-invariant because $\alpha$ has degree 0 and nowhere vanishing because $\Delta$ is. The modular unit $u_{\alpha}$ can be viewed as a morphism from the open modular curve of level $\Gamma_0(N)$ to the multiplicative group. 

\begin{definition} Given a divisor $\alpha$ as above, \emph{elliptic units} are values of $u_{\alpha}$ at points in $\Gamma_0(N) \backslash \mathcal H_{\infty}^{\mathrm{CM}}$. 
\end{definition}

Theorem \ref{T: main CM} implies that elliptic units are defined over ring class fields of imaginary quadratic fields. However, a  careful study of the integrality of modular units (cfr. \cite[Chp. 1, Lemma 1.1]{KL1981}) yields the following result.

\begin{theorem}
Let $\oo$ be an order in an imaginary quadratic field $K$, of discriminant prime to $N$. 
Let $u_{\alpha}$ be a modular unit of level $\Gamma_0(N)$ and let $\tau \in \Gamma_0(N)\backslash \cH_{\infty}^{\oo}$. Then 
\[u_{\alpha}(\tau) \in \mathcal O_H[1/N]^{\times},\] 
where  $\mathcal O_{H}$ is the ring of integers of the ring class field $H$. \end{theorem}

\begin{remark}In fact, the values $u_{\alpha}(\tau)$ for $\tau \in \Gamma_0(N)\backslash \cH_{\infty}^{\oo}$ can be used to produce genuine units in the ring class field of $\oo$. More precisely, one can show that for every $\sigma \in \mathrm{Gal}(H/K)$, we have
\[
(\sigma-1) u_{\alpha}(\tau) \in \cO_{H}^{\times}.
\] 
\end{remark}

\subsubsection{Heegner points}\label{SSS:Heegner points}

The construction of Heegner points on an elliptic curve $E/\Q$ relies crucially on the \emph{modularity} of the elliptic curve $E$. We briefly recall this notion.

Let $f=\sum_{n \geq 1} a_n q^n$ be a newform of weight 2 and level $\Gamma_0(N)$ for some $N\geq 1$; suppose that its Fourier coefficients are integers. The Eichler--Shimura construction attaches to $f$ an elliptic curve, as we now recall. Let $\mathbb T_2(\Gamma_0(N))$ be the Hecke algebra acting on cusp forms of weight $2$ and level $\Gamma_0(N)$ over $\Z$. Let $I_f$ be the ideal given by the kernel of the ring morphism 
 $\mathbb T_2(\Gamma_0(N)) \to \Z$ sending the Hecke operator $T_{n}$ to $a_n$. 
Let $J_0(N)$ be the Jacobian of the modular curve $X_0(N)$ of level $\Gamma_0(N)$. It is an abelian variety of dimension equal to the genus of $X_0(N)$, endowed with an action of $\mathbb T_2(\Gamma_0(N))$. The quotient 
\[
E_f:=J_0(N)/I_f J_0(N)
\]
is the elliptic curve attached to $f$. 
The composition of the morphism $X_0(N) \to J_0(N)$ given by $P \mapsto (P)-(\infty)$ with the projection $J_0(N) \to E_f$ yields a surjective $\Q$-rational morphism 
\begin{equation*} 
\pi_{E_f} \colon X_0(N) \to E_f. 
\end{equation*}	

\begin{definition} Given an elliptic curve $E$, we say that $E$ admits a modular parametrisation if there exists a surjective $\Q$-rational morphism $\pi_E \colon X_0(N) \to E$. When this is the case, we say that $E$ is \emph{modular}.
\end{definition} 
The work of Wiles \cite{Wiles}, Taylor--Wiles \cite{TW95} on Fermat's Last Theorem and its generalisations by Breuil, Conrad, Diamond, and Taylor \cite{BCDT} culminated in a full proof of the following theorem, formerly known as the Shimura--Taniyama Conjecture.
\begin{theorem}
Let $E/\Q$ be an elliptic curve. Then $E$ is modular. 
\end{theorem}
\begin{remark}
For an alternative formulation of modularity via $L$-functions see \cite{DDT1}.
\end{remark}

\begin{definition}
Let $E$ be a modular elliptic curve, and let $\pi_E \colon X_0(N) \to E$ be its modular parametrisation. \emph{Heegner points} on $E$ are images of points in $\Gamma_0(N) \backslash \cH_{\infty}^{\mathrm{CM}}$ under $\pi_E$.
\end{definition}

Theorem \ref{T: main CM} implies that Heegner points are defined over ring class fields of imaginary quadratic fields. Their importance in relation to the Birch and Swinnerton--Dyer conjecture will be discussed in \S \ref{SS: GZ Inventiones}.

	\section{Heegner constructions for real quadratic fields} \label{S: Heegner RM}
	Given the importance of singular moduli and more general Heegner constructions for algebraic number theory, it is desirable to have an analogue for real quadratic fields. 
	While no direct analogue is available for real quadratic fields, a (conjectural) theory of ``real multiplication'' relying on $p$-adic methods was recently proposed by Darmon and Vonk in \cite{DV2021}. In this section we will present the theory of \emph{rigid meromorphic cocycles}, which is the main new tool in this approach. A rigid meromorphic cocycle is a class in the first cohomology of the Ihara group $\mathrm {SL}_2(\mathbb Z[1/p])$ acting on the non-zero rigid meromorphic functions on the Drinfeld $p$-adic upper half plane $\mc{H}_p$.  
	These rigid meromorphic cocycles can be evaluated at real quadratic points of $\mc{H}_p$. Conjecturally their values display striking similarities with the classical singular moduli.  In a similar vein, the theory of analytic theta cocycles, obtained by considering cocycles for the Ihara group with values in analytic functions on $\cH_p$ defined up to scalars, can be used to reinterpret constructions of Darmon--Dasgupta \cite{DD06} and Darmon \cite{Dar01} of conjectural analogues of elliptic units and Heegner points for real quadratic fields. 
	
	Before introducing rigid meromorphic cocycles, we will review the relevant background on the Drinfeld $p$-adic upper half plane. 
	
	\subsection{Drinfeld $p$-adic upper half plane}\label{SS: Drinfeld upper half plane}
	In search for an analogue of singular moduli for real quadratic fields, the first obstacle one runs into is the absence of real quadratic points on the complex upper half plane. However, there is a $ p $-adic analogue of the complex upper half plane containing many real quadratic points, the \emph{Drinfeld upper half plane}, which we briefly describe. 
	The Drinfeld upper half plane $\mc{H}_p$ is a rigid analytic space\footnote{For an introduction to rigid analytic geometry see e.g.\ \cite{Conrad2008, bosch}.}, whose $ \C_p $-points are given by \[ \Hc_p=\mb{P}^1(\C_p) \setminus \mb{P}^1(\Q_p).\]
In the following we abuse notation and write $\mc{H}_p$ for the rigid space and as well as its $\C_p$-points. The space  $\mc{H}_p$ has an admissible covering by an increasing sequence of open affinoid subdomains $ \Hc_p^{\le n} $ which are constructed by removing balls of decreasing radius around all $ \Q_p $-rational points. More explicitly, for each $ z\in \mbP^1(\C_p) $ we fix homogeneous coordinates $ z=[z_0:z_1] $ such that  $\max\{ \vert z_0\vert, \vert z_1\vert \}=1 $. Then  \[ \Hc_p^{\le n}=\{ [z_0:z_1]\in \mb{P}^1(\C_p) \, \vert \ \order_p(x_0z_1-x_1z_0)\le n\ \forall [x_0:x_1]\in \mb{P}^1(\Q_p) \}. \]
	
For example, when $ n=0 $, the open affinoid $ \Hc_p^{\le 0} $ consists of the points $[z_0:z_1]\in \mbP^1(\C_p) $ such that $ \order_p(x_0z_1-x_1z_0)\le 0 $ for any $ [x_0:x_1]\in \mb{P}^1(\Q_p)  $. Hence  $ \Hc_p^{\le 0} $ is obtained from $\mb{P}^1(\C_p)$ by removing balls of radius~$1$ around the $ \Q_p $-rational points. The latter are precisely the points that are mapped to the $ \F_p $-rational points under the projection map $ \pi:\mbP^1(\C_p)\rightarrow \mbP^1(\overline{\mathbb{F}}_p) $, so that 
\[ \Hc_p^{\le 0}=\mbP^1(\C_p) \ \backslash \ \pi^{-1}(\mbP^1(\F_p)). \]

For convenience let us include the following \emph{ad hoc} description of the rigid analytic functions on $\mc{H}_p$, which the reader who is unfamiliar with rigid analytic geometry can take as a definition.
\begin{definition}\begin{itemize} 			

\item A \textit{rigid analytic function} on $ \Hc_p^{\le n} $ is a limit with respect to uniform convergence of rational functions on $ \mbP^1(\C_p) $ with poles in $ \mbP^1(\C_p) \ \backslash \ \Hc_p^{\le n} $.

\item A \textit{rigid analytic function} on $ \Hc_p $ is a function $ f:\Hc_p\rightarrow \C_p $, such that the restriction of $ f $ to $ \Hc_p^{\le n} $ is a rigid analytic function for every $ n\ge 0 $.
			
\item A \textit{rigid meromorphic function} on $ \Hc_p$ is a ratio $ f=\frac{g}{h} $ of rigid analytic functions $ g $ and $ h $ with $ h\neq 0 $.
\end{itemize}
	\end{definition}

	We denote by $ \mathcal{A} $ the set of rigid analytic functions on $ \Hc_p $ and by $ \mathcal{M} $ the set of rigid meromorphic functions.

The key structural properties of the rigid space $\mc{H}_p$ are that it is a smooth $1$-dimensional Stein space and that its algebra of rigid analytic functions $\mc{A}$ is a reflexive Fr{\'e}chet space. 

The group $\GL_2(\Q_p)$ acts on $\Hc_p $ by Möbius transformations, i.e., for $\gamma = \begin{psmallmatrix}
	a       & b \\
	c       & d
	\end{psmallmatrix} \in \GL_2(\Q_p)$ and $z \in \mc{H}_p$, the action of $\gamma$ on $z$ is given by 
	\[\gamma\cdot z = \frac{az+b}{cz+d}.\] 
As the center of $\GL_2(\Q_p)$ acts trivially, this also defines an action of $\PGL_2(\Q_p)$.

	\subsubsection{Bruhat--Tits tree for $ \PGL_2(\Q_p) $}
	
	A useful tool for studying the Drinfeld upper half plane is the Bruhat--Tits tree $\mc{T}$ attached to the group $ \PGL_2(\Q_p) $. We recall its construction.
	
	We call two $ \Z_p $-lattices $ \Lambda $ and $ \Lambda' $ in $ \Q_p^2 $ equivalent if they are homothetic. This defines an equivalence relation on the set of all $ \Z_p $-lattices in $ \Q_p^2 $. The set of vertices of the tree $ \mathcal{T} $ is defined to be the set of equivalence classes of lattices:\[ V(\mathcal{T}):=\{\Lambda\subseteq \Q_p^2\ \ \Z_p\text{-lattice} \}/\sim. \]
	By definition two vertices $ v,v'\in V(\mathcal{T}) $ are connected by an edge in the graph $ \mathcal{T} $ if there are representatives $ v=[\Lambda] $ and $ v'=[\Lambda'] $, such that $ p\Lambda\subsetneq\Lambda'\subsetneq \Lambda $.	One can show that the graph defined by this is a homogeneous tree of degree $ p+1 $. 
	For example, Figure 1  shows a part of the Bruhat--Tits tree for $ p=2 $.
	
	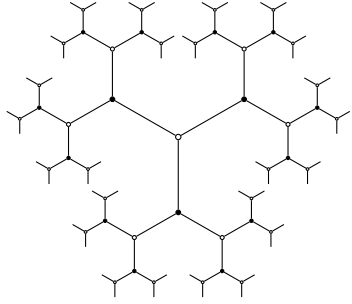
\begin{figure}
	\centering
		\begin{tikzpicture}[
		grow cyclic,
		level distance=1cm,
		level/.style={
			level distance/.expanded=\ifnum#1>1 \tikzleveldistance/1.5\else\tikzleveldistance\fi,
			nodes/.expanded={\ifodd#1 fill\else fill=none\fi}
		},
		level 1/.style={sibling angle=120},
		level 2/.style={sibling angle=120},
		level 3/.style={sibling angle=120},
		level 4/.style={sibling angle=120},
		nodes={circle,draw,inner sep=+0pt, minimum size=2pt},
		]
		\path[rotate=30]
		node []{}
		child foreach \cntI in {1,...,3} {
			node[circle,minimum size=1.75pt]{}
			child foreach \cntII in {1,...,2} { 
				node[circle,minimum size=1.5pt]{}
				child foreach \cntIII in {1,...,2} {
					node[circle,minimum size=1.25pt]{}
					child foreach \cntIV in {1,...,2} {
						node[circle,minimum size=1pt]{}
						child foreach \cntV in {1,...,2} {}
					}
				}
			}
		};
		\end{tikzpicture}
	\caption{The Bruhat--Tits tree of $\PGL_2(\Q_2)$ }
\end{figure}
	The group $ \GL_2(\Q_p) $ acts transitively on the set of all $ \Z_p $-lattices in $ \Q_p^2 $, preserving the equivalence classes. This induces a transitive action of $ \GL_2(\Q_p) $ and of $ \PGL_2(\Q_p) $ on $ \mathcal{T} $.
	
	\begin{example}
		The lattice $ \Z_p^2 $ gives rise to a vertex called the \textit{standard vertex} $  v_0=[\Z_p^2] $. The adjacent vertices of $ v_0 $ are given by \begin{align*}
		\begin{pmatrix}
		0	&-1  \\
		1	&0
		\end{pmatrix}
		\begin{pmatrix}
		1		&0	\\
		0		&p
		\end{pmatrix} v_0\, \text{  and  }\,\  
		\begin{pmatrix}
		1		&0  \\
		\lambda	&1
		\end{pmatrix} 	
		\begin{pmatrix}
		1		&0	\\
		0		&p
		\end{pmatrix}v_0\text{ for }\lambda \in \F_p.
		\end{align*}
		Moreover, if we consider the action of $ \SL_2(\Z[1/p]) $ on  $\mathcal{T}$, the stabiliser of the standard vertex is the group $ \SL_2(\Z) $ and the stabiliser of the edge $ e_0 $ from $ v_0 $ to $ \left (\begin{smallmatrix}
		1		&0	\\
		0		&p
		\end{smallmatrix} \right) v_0 $ is the congruence subgroup $ \Gamma_0(p) $.
	\end{example}
	
	There is the so-called reduction map	
	\[ \red:\Hc_p\rightarrow |\mathcal{T}| \] 
	from the $p$-adic upper half plane to the geometric realisation of $\mc{T}$. We refer to \cite{DT2007} for details regarding its construction. Via this map, the tree $ \mathcal{T} $ serves as a nice combinatorial tool to study the Drinfeld upper half plane. Let us indicate a few features of the reduction map. It is equivariant for the action of $ \PGL_2(\Q_p) $. Moreover we can recover the affinoids $ \Hc_p^{\le n} $ in the admissible covering described above as inverse images \[  \Hc_p^{\le n}=\red^{-1}(|\mathcal{T}^{\le n}|), \]
	where $ \mathcal{T}^{\le n} $ is the subtree of $ \mathcal{T} $ consisting of all vertices and edges which have distance $ \le n $ to the standard vertex $ v_0=[\Z_p^2] $. Finally, the reduction map allows us to view $ \Hc_p $ as a tubular neighbourhood of the Bruhat--Tits tree. For example, the inverse image of the edge $ e_0 $ is the annulus \[ \red^{-1}(e_0)=\{ z\in \C_p\vert 1< \vert z\vert < p \} .\] 

\subsubsection{RM points on the Drinfeld $p$-adic upper half plane} \label{subsubsec:RMpoints}

We introduce some notation for real quadratic points on $\mathcal H_p$ which will later  be used  to define the evaluation of rigid meromorphic cocycles.
\par 

Let $\tau \in \C_p$. We say that $\tau$ is a real quadratic point of discriminant $D$ if $[\tau: 1]$ satisfies the homogeneous quadratic equation $AX^2 + BXY+CY^2= 0$ with $A, B, C  \in \Z$ coprime integers satisfying $D = B^2-4AC > 0$, with $D$ not a square.

Then $K=\Q(\tau)$ is a real quadratic extension of $\Q$ and the point $[\tau:1] \in \mb{P}^1(\C_p)$ lies in $ \Hc_p $ if and only if the prime $p$ does not split in $K$. 
We refer to points  arising in this way as \emph{real multiplication} points or RM points of $\mc{H}_p$ and denote the set of these points as $\mc{H}_p^{\mathrm{RM}}$. 
For $\tau \in \mathcal{H}_p^{\RM}$, define 
	$$\mc{O}_\tau := \left\{ \begin{pmatrix}
	a       & b \\
	c       & d
	\end{pmatrix} \in \mathrm{M}_2(\Z [1/p]) : a \tau +b = c \tau^2 +d \tau \right\}. $$
	Then $\mc{O}_{\tau}$ is isomorphic to a $\Z [1/p]$ order in the real quadratic field $K=\Q(\tau)$ via the inclusion 
	$$\mc{O}_\tau \hookrightarrow K, \quad
	\gamma = \begin{pmatrix}
	a       & b \\
	c       & d
	\end{pmatrix}
	\mapsto c\tau +d.$$
Let $\Gamma:= \SL_2(\Z[1/p])$ denote the Ihara group. Let  $ \mc{O}^{\times}_{\tau,1}$ denote  the elements in $\mc{O}^{\times}_{\tau}$ of reduced norm one. The Dirichlet $S$-unit Theorem implies that the stabiliser 
\[\mathrm{Stab}_\Gamma(\tau) \cong \mc{O}_{\tau,1}^\times \cong \{\pm 1\} \times \langle \gamma_\tau \rangle\]
 is up to torsion a free group of rank $1$. It is generated by a fundamental unit $\gamma_{\tau}$ of norm $1$, called the \emph{automorph} of $\tau$. Moreover the points of $\mathcal H_p^{\RM}$ are precisely those for which the stabiliser in $\Gamma$ has rank 1. This feature is crucial to the evaluation of the rigid meromorphic cocycles defined in the subsequent sections. We note that the action of $\Gamma$ on $\mathcal H_p^{\RM}$ need not preserve the discriminant, but  it does preserve its prime-to-$p$ part.

\subsection{Group actions, curves and $ p $-adic theta functions}\label{subsection 3.2}

We review the theory of $p$-adic theta functions arising in the context of $p$-adic uniformisation of Shimura curves. While the construction of rigid cocycles does not fit in this setting, the cocycles arising in \S \ref {SS:  RMC} 
are obtained via similar methods, after an appropriate degree shift in cohomology. \\
\par

Let $ \Gamma_0\subset \SL_2(\Q_p)$ be any subgroup. We say that $\Gamma_0$ acts discretely on $\mc{H}_p$ if the orbit $ \Gamma_0 \tau $ of each element $ \tau\in \Hc_p $ intersects each of the affinoids $ \Hc_p^{\le n} $ only in finitely many points. 

		Note that the groups $ \SL_2(\Z) $, $ \SL_2(\Z[1/p]) $ and $ \SL_2(\Z_p) $ all act non-discretely on $ \Hc_p $. To see this recall that $ \SL_2(\Z) $ is contained in the stabiliser of the standard vertex $ v_0 $ of the Bruhat--Tits tree and that any element $ \tau\in\Hc_p^{\le 0} $ is mapped to $v_0$ via the reduction map. So the whole orbit $ \SL_2(\Z)\tau $ is contained in $ \Hc_p^{\le 0}=\red^{-1}(\mathcal{T}^{\le 0}) $ and the set $ \SL_2(\Z)\tau\cap\Hc_p^{\le0} $ contains infinitely many elements. 

The upper half plane $\mc{H}_p$ is of major significance for arithmetic geometry as it can be used to construct families of algebraic curves (so-called Mumford curves) which uniformise certain Shimura curves $p$-adically. 
These curves arise as quotients of $\mc{H}_p$ by discrete group actions. We refer to \cite[\S 2.4]{DT2007} and references therein for details. Here we only mention the following class of examples. 
	
	\begin{example}\label{example of Gamma acting discretely on H_p}
		Let $ B $ be the definite quaternion algebra $ \Q[i,j] \subset \mathbb{H}$ inside the Hamilton quaternions and let $ R=\Z[i,j,\frac{i+j+k+1}{2}] $ be a maximal order. Choose an odd prime $ p $ at which $ B $ splits and consider the group $ \Gamma_B=(R[1/p])_1^\times $ of units $ \gamma \in R[1/p]^\times$ of reduced norm $1$. After choosing an isomorphism $(B\otimes_{\Q}\Q_p)^{\times} \cong \GL_2(\Q_p)$ the group $ \Gamma_B $ embeds into $\SL_2(\Q_p)$ and acts discretely on $ \Hc_p $. Let $ \Gamma'\subset \Gamma_B $ be a subgroup of finite index which has no fixed points on $ \Hc_p $. Then the quotient space $ X_{\Gamma'}=\Gamma'\backslash\Hc_p $ is a rigid analytic space which admits a model as a Shimura curve.
	\end{example}
	
	Note that a definite quaternion algebra does not admit an embedding of a real quadratic field. Essentially due to this fact the quotients $ X_{\Gamma'}$ are not so relevant for a theory of real multiplication. As we will see below, a better setup for these purposes is to consider the indefinite quaternion algebra $ B=\mathrm{M}_2(\Q) $ and the (non-discrete) actions of groups like $ \SL_2(\Z) $ and $ \SL_2(\Z[1/p]) $ on $\mc{H}_p$. 
	
Nevertheless we want to briefly discuss one aspect from the theory of Mumford curves, the construction of $p$-adic theta functions. These are important tools in the study of the Jacobian of the compactification of a Mumford curve. Key ideas of their construction can be modified to construct interesting functions also in situations where a non-discrete group action is involved. 
Assume that $ \Gamma_0\subset \SL_2(\Q_p)$ is a subgroup which acts discretely and without fixed points on $\mc{H}_p$. The $p$-adic theta functions for $\Gamma_0$ are meromorphic functions on $ \Hc_p $, which are invariant under $ \Gamma_0 $ up to multiplicative scalars and which can be constructed as follows.
	
For an element $ w\in\mb{P}^1(\C_p) $, define the rational functions $ t_w$ on $\mb{P}^1(\C_p)$ by \[ t_w(z)=\begin{cases}
	z-w, \ \ \ \ \  \text{ if }\ \vert w \vert \le 1,\\
	z/w-1,\ \text{ if }\ \vert w\vert>1,\\
	1, \ \ \ \ \ \ \ \ \ \ \  \text{ if }w=\infty.
	\end{cases} \] 
	Note that for two distinct elements $ w^+ ,w^-$ in $ \Hc_p $ and $\gamma \in \Gamma_0$ the quotient function \[ \frac{t_{w^+}(\gamma z)}{t_{w^-}(\gamma z)} \] has a simple zero at $ \gamma^{-1}w^+ $ and a simple pole at $ \gamma^{-1}w^- $, as does the function $ \frac{t_{\gamma^{-1}w^+}( z)}{t_{\gamma^{-1}w^-}( z)}$. Therefore they are equal up to a constant. 
	
	Now consider the degree zero divisor $ \Delta=w^+-w^- $. Then the product 

	\[ f_\Delta(z)=\prod_{\gamma\in \Gamma_0} \frac{t_{\gamma w^+}(z)}{t_{\gamma w^-}(z)} \]
	converges to a rigid meromorphic function on $ \Hc_p $. (For more details, see Section 2.2 in \cite{DV2021}.) By the observation above, the function $ f_\Delta $ is $ \Gamma_0 $-invariant up to multiplication by a constant in $ \C_p $. Hence, it defines an element of $ (\mathcal{M}^\times/\C_p^\times)^{\Gamma_0} =H^0(\Gamma_0,\mathcal{M}^\times/\C_p^\times)$. 
The function $f_\Delta$ is called the \emph{$p$-adic theta function} associated to the degree zero divisor $\Delta$.

	\subsection{Rigid meromorphic cocycles}\label{SS:  RMC} 
	
	One way of looking at the $j$-function is as an $\SL_2(\Z)$-invariant function on the complex upper half plane $\mc{H}_{\infty}$, i.e., an element of $H^0( \SL_2(\Z), \mathrm{Hol}(\mc{H}_{\infty})),$
	where  $\mathrm{Hol}(\mc{H}_{\infty})$ denotes the holomorphic functions on $\mc{H}_{\infty}$.
	Now consider $ \Gamma:=\SL_2(\Z[1/p]) $ acting on $\mc{H}_p$. In search for direct analogues of $j$ in the $p$-adic setting, one might first naturally want to consider the groups $H^0(\Gamma,\mathcal{A})$ or $H^0(\Gamma,\mathcal{M})$. However it follows essentially from the Weierstrass preparation theorem that any $\Gamma$-invariant rigid analytic function must be constant (cf. \cite[Lemma 1.9]{DV2021}), hence $H^0(\Gamma,\mathcal{A})\cong \C_p \cong H^0(\Gamma,\mathcal{M})$ and so these groups do not contain interesting functions. Considering the corresponding cohomology groups in degree one on the other hand turns out to be a successful strategy, even more so if one furthermore switches to the multiplicative theory. This leads to the following definition:
	\begin{definition}
		\begin{itemize}
			\item A \textit{rigid meromorphic} (resp.~\textit{analytic}) \textit{cocycle} is an element in $ H^1(\Gamma,\mathcal{M}^\times) $ (resp. in $ H^1(\Gamma,\mathcal{A}^\times) $).
			\item A \textit{rigid meromorphic} (resp.~\textit{analytic}) \textit{theta cocycle} is an element in $ H^1(\Gamma,\mathcal{M}^\times/\C_p^\times) $ (resp.~in $ H^1(\Gamma,\mathcal{A}^\times/\C_p^\times) $).
		\end{itemize}
	\end{definition}
Let $\Gamma_{\infty}\subset \Gamma$ be the subgroup of upper triangular matrices. A rigid meromorphic cocycle is called \emph{parabolic} if its restriction to $\Gamma_{\infty}$ is trivial. It is called \emph{quasi-parabolic} if its restriction to $\Gamma_{\infty}$ lies in the image of $H^1(\Gamma_{\infty}, \C_p^{\times})$. The groups of such cohomology classes are denoted by $H^1_{\mathrm{par}}(\Gamma,\mathcal{M}^\times)$ and $H^1_f(\Gamma, \mc{M}^{\times})$, respectively. One can show that any cohomology class in $H^1_f(\Gamma, \mc{M}^{\times})$ has a unique representative, whose values on $\Gamma_{\infty}$ consist of constant functions.

	\begin{remark}	 One reason for considering the multiplicative group $ \mathcal{M}^\times $ as coefficients rather than the additive group $ \mathcal{M} $ is the connection to the theory of modular symbols and so-called rigid meromorphic period functions. More precisely, the subgroup $ H^1_{\mathrm{par}}(\Gamma,\mathcal{M}^\times)$ can be identified with the group of $\Gamma$-invariant modular symbols $\mc{MS}^{\Gamma}(\mc{M}^{\times})$ (\cite[Section 1]{DV2021}). This connection to modular symbols is useful for establishing structural results. It also allows one to produce many interesting elements in $ H^1(\Gamma,\mathcal{M}^\times) $. Contrary to this, the group $ H^1(\Gamma, \mathcal{M}) $ contains no \textit{parabolic} elements (\cite[Remark 1.8]{DV2021}). 
	\end{remark}	
	
In practice it is easier to write down examples of theta cocycles (we will give a few of these below). Note that there is a natural map $ H^1(\Gamma,\mathcal{M}^\times) \rightarrow  H^1(\Gamma,\mathcal{M}^\times/\C_p^\times) $ and one would like to lift  theta cocycles to rigid meromorphic cocycles. Unfortunately this is in general not possible as there is an obstruction to lifting. What is possible though is to lift the restriction of a theta cocycle to the group $\SL_2(\Z)$ to an element in $H^1(\SL_2(\Z),\mathcal{M}^\times)$, as $H^2(\SL_2(\Z),\C_p^\times)=0$. In some situations this is good enough. 
\newpage
	
	\subsubsection{Examples}
	Let us discuss some examples. 
	\begin{example}[The trivial cocycle]
		Fix a base point $ \xi\in \mbP^1(\Q_p) $. For any $ \gamma \in \Gamma $, define the rational function 
		\[ J_{\triv}(\gamma)(z):=\frac{z-\gamma\xi}{z-\xi}. \] 
		
Up to a scalar this function is uniquely determined as the rational function with divisor $ (\gamma\xi)-(\xi) $. Then for two elements $ \gamma_1 $, $ \gamma_2\in \Gamma $, the function $ J_{\triv}(\gamma_1\gamma_2) $ is the function with divisor $ (\gamma_1\gamma_2\xi)-(\xi) =(\gamma_1\gamma_2\xi)-(\gamma_1\xi)+(\gamma_1\xi)-(\xi)$. Hence up to a scalar, $ J_{\triv}(\gamma_1\gamma_2) $ is the same as $ J_{\triv}(\gamma_1) \gamma_1J_{\triv}(\gamma_2)  $. So we obtain a rigid analytic theta cocycle $J_{\triv}\in H^1(\Gamma, \mathcal{A}^\times/\C_p^\times) $. 
	\end{example}

	\begin{example} (Theta cocycles attached to RM points). This example is inspired by a similar construction in the theory of rational modular cocycles (see \cite[Section 1.4]{DV2021}).
		As we have observed above, $ \Gamma $ does not act discretely on $ \Hc_p $. However, it does act discretely on the product $ \Hc_p\times \Hc_\infty $ of the Drinfeld upper half plane and the complex upper half plane $ \Hc_\infty $. We can use this to construct examples of theta cocycles as follows.
		
		For a tuple $ (r,s)\in\mbP^1(\Q)^2 $ and $ w $ real quadratic, let $ w' $ be its algebraic conjugate. The geodesic from $ r $ to $ s $ in $ \Hc_\infty $ cuts ${\overline{\Hc}}_\infty=\Hc_\infty\cup\R\cup\{\infty\} $ into two connected components, say "on the right hand side" of the geodesic is the $ - $-component $ \Hc_\infty^- $ and on the left hand side is the $ + $-component $ \Hc_\infty^+ $ (see Figure \ref{intgeonumber}). This allows us to define the following symbol
		\begin{equation}\label{intgeo}
			 (r,s)(w,w'):=\begin{cases}
		0\text{, if }w,w' \text{ are in the same connected component}\\
		1\text{, if }w\in \Hc_\infty^+,\ w'\in \Hc_\infty^- \\
		-1\text{, if }w\in \Hc_\infty^- ,\ w'\in\Hc_\infty^+.
		\end{cases}
		\end{equation}	
		We refer to it as the signed intersection number of the geodesics from $ r $ to $ s $ and from $ w $ to $ w' $. 
		For example the intersection number of the geodesics in Figure \ref{intgeonumber}  is $ -1 $.
		\begin{figure}[h!] 
		\begin{equation*}
		\begin{tikzpicture}
		\draw (0,0) -- (8,0);
		\draw (5,0) arc (0:180:2cm);
		\draw (7,0) arc (0:180:1.5cm);
		\node[below] at (1,0) {$ r $};
		\node[below] at (5,0) {$ s $};
		\node[below] at (4,0) {$ w $};
		\node[below] at (7,0) {$ w' $};
		\node[below] at (3,1) {$ - $};
		\node[below] at (3,3) {$ + $};
		\node[above] at (7.7,0) {$ \Hc_\infty $};
		\end{tikzpicture}
		\end{equation*}
\caption{Intersection of geodesics}
\label{intgeonumber}
		\end{figure}
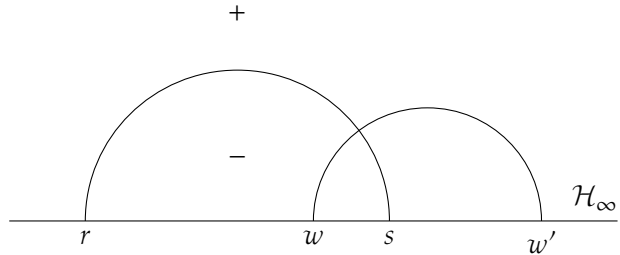

			For an RM point $ \tau \in \Hc_p^{RM} $,  define the function 
			\begin{align*}
			J_\tau\colon & \Gamma\rightarrow\mathcal{M}^\times\\
			&\gamma\mapsto \prod_{w\in\Gamma\tau}t_w(z)^{(\gamma\infty,\infty)(w,w')},
			\end{align*}
where $t_w(z)$ is as in Section \ref{subsection 3.2}. 

\noindent As it turns out, the set of elements $ w \in\Gamma\tau$ such that 
$ (\gamma\infty,\infty)(w,w') \neq 0$
 is a discrete subset of $ \Hc_p $. One can show that the product converges, that indeed $ J_\tau(\gamma) $ defines a non-zero rigid meromorphic function and that furthermore $ J_\tau $ defines a class in $ H^1(\Gamma,\mathcal{M}^\times/\C_p^\times) $.
	\end{example}
	\begin{example}[The winding cocycle]
		Let us also mention here the so-called \textit{winding cocycle} from \cite[Section 2.3]{DPV1}, which is an example of a rigid analytic theta cocycle. 
		To define it choose a base point $ \xi = (\xi_p, \xi_\infty)\in \Hc_p\times \Hc_\infty $, such that $ \xi_p\in \Hc_p^{\le 0} $ and $ \xi_\infty $ does not lie in any geodesic in the $ \Gamma $-orbit of $  [0,\infty]$. Denote by $ \Sigma $ the $ \Gamma  $-orbit of the pair $ (0,\infty) $ in $ \mb{P}^1(\Q)^2$ and define for every $ \gamma \in \Gamma $ and $z\in \mc{H}_p$ 
		\[ J^\xi_{(0,\infty)}(\gamma)(z):=\prod_{(r,s)\in\Sigma}\left(\frac{\xi_p-r}{\xi_p-s} \cdot \frac{z-s}{z-r}\right)^{(r,s)(\xi_\infty,\gamma \xi_\infty)}, \]
		using the usual conventions when some of the points are $\infty$ and 
		where $ (r,s)(\xi_\infty,\gamma \xi_\infty) $ denotes the signed intersection number as in (\ref{intgeo}). 
		The function $ J^\xi_{(0,\infty)}(\gamma) $ converges to a rigid analytic function on $ \Hc_p $ and $\gamma \mapsto  J^\xi_{(0,\infty)}(\gamma)$ satisfies the cocycle condition up to scalars in $ \C_p^\times $ (cf.~\cite[Proposition 2.5]{DPV1}). The corresponding  class in $ H^1(\Gamma, \mathcal{A}^\times/\C_p^\times) $ is independent of the choice of base point \cite[Prop. 2.6]{DPV1} and is called the winding cocycle and we denote it by $J_{(0,\infty)}$.	
	\end{example}

	\subsubsection{Hecke module structure of rigid meromorphic cocycles}
	On the space $H^1(\Gamma, \mathcal{M}^\times)$ of rigid meromorphic cocycles and on the space $H^1(\Gamma, \mathcal{M}^\times/ \mathbb{C}_p^\times)$ of rigid meromorphic theta cocycles one can define Hecke operators $T_n$ for any $n\geq 1$. This is done by decomposing a relevant union of double cosets into cosets and following a recipe of Shimura. We refer to \cite[\S 2.4]{DPV1} for details. 

\begin{definition}\label{D: period function}
Given a quasi-parabolic rigid meromorphic cocycle $J$, its  \emph{rigid meromorphic period function} is the value at $S:=\begin{psmallmatrix} 0 & -1 \\ 1 & 0 \end{psmallmatrix} \in \Gamma$ of the quasi-parabolic representative of $J$.
\end{definition}

Let us define the abstract Hecke algebra as $\mathbb{T} := \Z[T_n \ , \ n\geq 1]$. One can then prove the following structural result.
	\begin{theorem}[{\cite[Theorem 1 and Remark 2]{DV2021}}]  The group
		$H^1(\Gamma, \mathcal{M}^\times) $ is of infinite rank. Any finite rank  $\mathbb{T}$-stable submodule of $H^1_f(\Gamma, \mathcal{M}^\times)$ is contained in $H^1_f(\Gamma, \mathcal{A}^{\times})$.
	\end{theorem}
	 
Analytic theta cocycles behave differently. As Hecke modules, they are closely related to modular forms. In \cite{DV22}, Darmon and Vonk prove a classification result using ideas of Stevens and Schneider--Teitelbaum. We briefly review the key points (see also \cite[Section 3.1]{DPV1}).\\
	\par
  	Consider the annulus $U := \{ z \in \mathbb{C}_p| \ 1 < |z| < p\}$. Its stabiliser in $\Gamma$ is $\Gamma_0(p)$. Denote $\text{dlog}f: = \frac{f'(z)}{f(z)}dz.$ The group homomorphism 
    	\[
    	\mathcal A^{\times} \to \mathbb C_p, \qquad f \mapsto {\res_{U}}(\text{dlog} f), 
    	\]
 where $\res_{U}$ denotes the $p$-adic annular residue along $U$, is $\Gamma_0(p)$-equivariant and trivial on constants. It can be shown that it takes values in $\Z$. This induces a map on cohomology 	
	\[
\delta_U:=	\res_U\circ \text{dlog} \colon 
H^1(\Gamma, \mathcal{A}^\times/ \mathbb{C}_p^\times) \rightarrow H^1(\Gamma_0(p),\Z).
	\]	
	
	\begin{theorem}\cite[Section 3]{DV22}.
	The map 
	\[
	\delta_U \otimes 1: H^1(\Gamma, \mathcal{A}^\times/ \mathbb{C}_p^\times)\otimes \Q \rightarrow H^1(\Gamma_0(p),\Z)\otimes \Q
	\]
			is surjective and has a $1$-dimensional kernel, generated by the analytic theta cocycle $J_{\triv}$. 
	
	The map $\delta_U$ admits a Hecke-equivariant section
			$$\ST:H^1(\Gamma_0(p),\Z) \rightarrow H^1(\Gamma, \mathcal{A}^\times/ \mathbb{C}_p^\times).$$
	\end{theorem}
	
The section $\ST$ above is called the multiplicative \emph{Schneider--Teitelbaum lift}, and one can use it to classify the analytic theta cocycles as follows.

		\begin{theorem}\cite[Section 3]{DV22}. The space
	$H^1(\Gamma, \mathcal{A}^\times/ \mathbb{C}_p^\times) \otimes \Q$ is of dimension $2g+2$, where $g$ is the genus of the modular curve $X_0(p)$. The Hecke action factors through the faithful Hecke algebra of modular forms $M_2(\Gamma_0(p))$.
	\end{theorem}

\subsection{Heegner constructions} \label{SS: Heegner objects RM}

The crucial motivation for the study of rigid meromorphic cocycles and analytic theta cocycles is the construction of conjectural analogues of Heegner objects for real quadratic fields. In analogy with the Heegner constructions, they can be viewed as images of maps 
\[
\Phi \colon \Gamma \backslash \cH_p^{\mathrm{RM}} \to X(\C_p) 
\]
for the group schemes $X/\Q$ appearing in \S \ref{SS: Heegner Objects CM}.

 The corresponding images are expected to be global. Note that however in the RM settings these maps will be defined exclusively at RM points in $\mathcal H_p$ via the evaluation of rigid meromorphic or theta cocycles, which we now describe.
\subsubsection{Evaluation at RM points} 
We now define the evaluation of quasi-parabolic rigid meromorphic cocycles at RM points. Recall that for $\tau \in \mathcal{H}_p^{\RM}$, $\gamma_{\tau} \in \Gamma$ denotes the automorph of $\tau$ as defined in \ref{subsubsec:RMpoints}. 
	
	\begin{definition}
		Let $[J] \in H^1_f(\Gamma, \mathcal{M}^\times)$ be a quasi-parabolic rigid meromorphic cocycle. Consider a representative $J \in [J]$ given by a map  $J\colon \Gamma \rightarrow\mathcal{M}^\times$ such that $J(\Gamma_{\infty})\subset \C_p^{\times}$ and let $\tau \in \mathcal{H}_p^{\RM}$. Define the \emph{evaluation of $J$ at $\tau$} as 
		$$ J[\tau] := J(\gamma_\tau)(\tau) \in \mathbb{C}_p \cup \{ \infty \}.$$
	\end{definition}
	
		\begin{remark} One can check that this is constant on $\Gamma$ orbits. 
	\end{remark}
	\begin{remark} Defining the evaluation of a theta cocycle at an RM point is also possible but requires some extra work involving the restriction to (a conjugate of) $\SL_2(\Z)$. We refer to Section 2.3 of \cite{DV22} for the definition. When the theta cocycle lifts to a rigid meromorphic cocycle, the values essentially agree. When the cocycle does not lift, the values are expected to be transcendental in general, but are still of arithmetic interest as they are related to Gross--Stark units (see Conj. \ref{C: Darmon-Dasgupta}) and Stark--Heegner points (see Conj. \ref{C: Stark-Heegner}) below. 
	\end{remark}

\subsubsection{Singular moduli for real quadratic fields}

Let $\tau \in \mathcal{H}_p^{\RM}$ and write $H_{\tau}$ for the  narrow ring class field corresponding to $\mc{O}_\tau$.  For a quasi-parabolic rigid meromorphic cocycle $J$ we let $j$ be its associated period function (see Def.\ \ref{D: period function}). Then we define a field
\[
H_j=H_J:=\text{Compositum}_{j(\tau)=\infty}(H_{\tau}),
\]
which is a finite extension of $\Q$.
\begin{definition}
Let $\tau \in \Gamma \backslash \cH_{p}^{\mathrm{RM}}$. For a quasiparabolic meromorphic cocycle $J$, we say that $J[\tau]$ is a \emph{real quadratic singular modulus}.
\end{definition}

 The following conjecture, due to Darmon and Vonk, justifies the analogy with the values of the $j$-function at CM points. 
	
\begin{conjecture}[\cite{DV2021}] If $J \in H^1_f(\Gamma, \mc{M}^{\times})$, then $J[\tau]$ is algebraic and lies in the compositum of $H_J$ and $H_{\tau}$. 
\end{conjecture}
Note that the analogy with classical singular moduli is imperfect, because the choice of cocycle $J$ impacts the field of definition of its value. A perhaps more appropriate analogy is that these values should be thought of as multiplicative analogues of differences of singular moduli, lying in composita of the corresponding ring class fields.	

\subsubsection{Elliptic units for real quadratic fields} 

Consider the unique (normalised) Eisenstein eigenform of weight 2 and level $\Gamma_0(p)$ with Fourier expansion
\[
E^{(p)}_2=\frac {p-1}{12}+\sum_{n \geq 1}\left (\sum_{d \mid n, p \nmid d} d \right  )q^n.  
\] 
For a choice of base point $z_0 \in \cH_{\infty}$, the function computing path integrals of the differential attached to the Eisenstein series 
\[ \gamma \mapsto \int_{z_0}^{\gamma z_0}\frac{1}{2\pi i}E^{(p)}_2(z)dz 
\]
for $\gamma \in \Gamma_0(p)$ takes integer values and transforms as a (non-parabolic) cocycle $\varphi_{\mathrm{DR}}  \in H^1(\Gamma_0(p), \Z)$. 
Its values can be described explicitly in terms of Dedekind sums (see \cite[\S 2.5]{DD06}).

The Schneider--Teitelbaum lift 
\begin{equation*} \label{E: J_DR definition}
J_{\mathrm{DR}}:={\ST} (\varphi_{DR})
\end{equation*}
of $\varphi_{DR}$ to $H^1(\Gamma, \mathcal A^{\times}/\C_p^{\times})$
is called the \emph{Dedekind--Rademacher theta cocycle.}
\begin{definition}
Let $\tau \in  \Gamma \backslash \mathcal{H}_p ^{\RM}$ be a point of discriminant $D$ prime to $p$. The value $J_{\mathrm{DR}}[\tau]$ is called the \emph{real quadratic elliptic unit } attached to $\tau$.
\end{definition}

 Its values at RM points in $\mathcal H_p$ should be thought of as analogues of the elliptic units defined in \S \ref {S: Heegner CM}, as the following conjecture suggests.  

\begin{conjecture}\cite[Conj. 2.14]{DD06} 
\label{C: Darmon-Dasgupta}
 Let $\tau \in \cH_{p}$ be a real quadratic point of discriminant $D$ prime to $p$ and let $H_{\tau}$ be the narrow ring class field of the order defined by $\tau$. Then the real quadratic elliptic unit attached to $\tau$ satisfies 
\[
J_{\mathrm{DR}}[\tau] \in \cO_{H_{\tau}}[1/p]^{\times}.
\]	
\end{conjecture}

Note that unlike in the CM setting, these values are expected to be $p$-units rather than genuine units in the Hilbert class field. The analogy with the construction of elliptic units is perhaps not obvious on the surface. However, the Eisenstein series $E_{2}^{(p)}$ whose periods appear in the definition above can be viewed as the logarithmic derivative of $\Delta(pz)/\Delta (z)$, a typical example of the modular unit featuring in \S \ref{SSS: elliptic units}.

\subsubsection{Stark--Heegner points}

Let $f$ be a newform of weight $2$ and level $\Gamma_0(p)$, with Fourier coefficients valued in $\Z$. Consider the differential form $\omega_f = 2 \pi i f(z) dz$, and let 
\[
\omega_f^{+} = \omega_f + \overline{\omega}_f,  \text{ and } \omega_f^{-} = \omega_f -\overline{\omega}_f.
\]
be its real and imaginary part, respectively. Given $z_0 \in \mathcal{H}_\infty$, the maps $\tilde {\phi}_f \colon \Gamma_0(p) \to \R$ defined by  
\[
\gamma \mapsto  \int^{\gamma z_0}_{z_0} \omega_f^\pm,
\]
can be rescaled by suitable periods $\Omega^\pm_f$ in such a way that the pair of cocycles $\varphi_f^{\pm}:=1/\Omega_{f}^{\pm}\tilde \varphi_f^{\pm} $ take values in $\Z$. The Schneider--Teitelbaum lifts of the cocycles defined above 
\[
J_f^{\pm}:=\ST(\varphi_f^{\pm})
\]	
are called the \emph{elliptic theta cocycles} attached to $f$. 
	
By the Eichler--Shimura construction reviewed in \S \ref{SSS:Heegner points}, a newform $f\in S_2(\Gamma_0(p))$ as above corresponds to an elliptic curve $E:=E_f$ with multiplicative reduction at $p$. For such an elliptic curve, Tate showed that its $\C_p$-valued points can be uniformised as
$$\Psi_\mathrm{Tate} \colon \mathbb{C}_p^\times/q_E^{\Z} \xrightarrow{\sim} E(\mathbb{C}_p)$$
for some parameter $q_E \in \mathbb{C}_p^\times$ satisfying $|q_E| < 1$. 

\begin{definition}
 Let  $\tau \in  \Gamma \backslash \mathcal{H}_p ^{\RM}$ be a point of discriminant $D$ prime to $p$. A \emph{Stark--Heegner point} for $E$ attached to $\tau$ is 
\[
P_\tau^{\pm}:=\Psi_{\mathrm{Tate}}(J_f^{\pm}[\tau]) \in E(\mathbb{C}_p).
\]
\end{definition}

The following conjecture, originally formulated by Darmon in \cite{Dar01} motivates the conjectural analogy between Stark--Heegner points and the classical Heegner points in CM theory. 

\begin{conjecture}[\cite{Dar01},{\cite[Conjecture 3.18]{DV22}}] \label{C: Stark-Heegner}
Let  $\tau \in \Gamma \backslash \mathcal{H}_p ^{\RM}$ be a point of discriminant $D$ prime to $p$.  Denote  $H_{\tau}$ the narrow ring class field of the order $\oo_\tau$ attached to $\tau$. Then 
\[P_{\tau}^{\pm} \in E(H_{\tau}). \]
\end{conjecture}				
In addition, the point $P_{\tau}^{\pm}$ is expected to lie in the $(\pm1)$-part of $E(H_{\tau})$ for the action of the complex conjugation in the narrow Hilbert class field of $\tau$. 

\begin{remark}
The above conjecture (and suitable generalisations) has been verified extensively, see for example \cite{Darmon-Pollack}. The theoretical evidence is much more fragmentary, and no general approach is known at this stage. However, an important case was showed in the article \cite{Bertolini-Darmon} for Stark--Heegner points defined over genus fields of real quadratic fields by relating them to Heegner points in this special setting. 
\end{remark}

\section{Families of real analytic Eisenstein series and CM theory} \label{S: analytic families}

	 In this section, we introduce the second main theme of this article, which is how modular generating series constructed from Heegner objects appear in relation to derivatives of families of modular forms. 
We will start by reviewing some archimedean examples; we will then move to their $p$-adic counterparts in \S \ref{S: p-adic families}. 

	In the archimedean setting, derivatives of real analytic families of Eisenstein series fit prominently in the Kudla program. Roughly speaking, this vast program seeks to show that certain formal $q$-expansions encoding special cycles on orthogonal Shimura varieties are modular. The most basic instances of these statements involve CM cycles on modular curves. The works \cite{GZ1985}, \cite{GZ1986} of Gross and Zagier were indispensable to the emergence of this program.  	
	 In \S \ref{SS: GZ Crelle}, we will give an account of their work on the differences of singular moduli, which will be close in spirit to the $p$-adic analogues discussed in \S \ref{S: p-adic families}. We will then give a brief overview of the work of Gross and Zagier on the derivative of the $L$-function attached to an elliptic curve \cite{GZ1986}. Finally in \S  \ref{SS: biquadratic framework} we explain how the two articles can be interpreted under a unified perspective, which will be translated in the $p$-adic setting in \S\ref{S: p-adic families}. 

\subsection{Gross--Zagier's factorisation formula} 
\label{SS: GZ Crelle}
We review the work of Gross and Zagier on the factorisation of differences of singular moduli, the CM values of the $j$-function introduced in \S \ref{SSS: singular moduli}.  Their landmark work unveiled new phenomena in the theory of complex multiplication, initiating the study of the intersection theory of CM cycles on modular curves via automorphic methods. 

	Let us introduce the main result. Let $D_1$ and $D_2$ be two negative relatively prime fundamental discriminants. Let  $D=D_1 D_2$ and let $w_i$ be the number of roots of unity in the quadratic orders of discriminant $D_i$.	
	For $\tau\in \mathcal{H}_{\infty}$ denote by $[\tau]$ the equivalence class of $\tau$ under the action of $\SL_2(\Z)$. The main theorem of complex multiplication implies that  the product 
	\begin{align*}
	\label{jproduct}
	J(D_1,D_2)=\prod_{\substack{[\tau_1],\, \,\disc \tau_1=D_1 \\ [\tau_2],\, \, \disc \tau_2=D_2}} \left(j(\tau_1)-j(\tau_2)\right)^{\frac{4}{w_1 w_2}} \ \ 
	\end{align*}
	(or more precisely, its square for $D_i\geq -4)$
	is an integer. Gross and Zagier are motivated by the deceptively simple question of determining its prime factorisation.
		
	For primes $\ell$ with $\left( \frac{D_1 D_2}{\ell} \right) \neq - 1$ define
	\[
	\epsilon(\ell)= \begin{cases}
	\left( \frac{D_1}{\ell} \right)  & \text{ if } (\ell,D_1)=1, \\
	\left( \frac{D_2}{\ell} \right)  & \text{ if } (\ell,D_2)=1.\\
	\end{cases}
	\]
	Now if $n=\prod_i \ell_i^{a_i}$ is such that $\left( \frac{D_1 D_2}{\ell_i} \right) \neq -1$ for all $i$, then we define $\epsilon(n)=\prod_i \epsilon(\ell_i)^{a_i}$.
	Gross and Zagier prove the following theorem.
	
	\begin{theorem}\cite[Thm 1.3]{GZ1985} 
		\begin{equation}\label{E: factorisation}
		J(D_1,D_2)^2= \pm \prod_{\substack{x,n,n' \in \Z \\ n,n'>0, \\ x^2+4nn'=D}} n^{\epsilon(n')}.
		\end{equation}
	\end{theorem}
	
	In \emph{loc. cit}, the authors provide two alternative proofs of this result. The first one is essentially algebraic and relies on  the moduli interpretation of CM points on the modular curve. The second is purely analytic and makes use of automorphic techniques to obtain a formula for $J(D_1, D_2)^2.$ While one can easily deduce the algebraicity of $J(D_1, D_2)$ from the theory of complex multiplication, no use of this fact is made in the analytic proof, where instead it can be seen as a byproduct of the proof itself. In view of translating the Gross--Zagier approach to the RM setting, the analytic strategy is of course more relevant, as a moduli interpretation of RM points of the Drinfeld $p$-adic upper half plane is lacking. \\

\par We will briefly outline some ideas behind the algebraic proof before focusing on the analytic one. The general idea is rather natural:  if a prime $p$ divides $J(D_1, D_2)^2$, there is a pair   of elliptic curves with CM by two orders $\mathcal{O}_{K_1}$ and $\mathcal{O}_{K_2}$ with isomorphic reductions over $\overline{\F}_p$.  This forces the endomorphism ring of the reduction to be a maximal order in the quaternion algebra ramified at $p$ and $\infty$. 

The existence of embeddings of the orders $\cO_{K_i}$ into an order in a quaternion algebra of discriminant $p$ is rather restrictive on the prime. Working out explicit conditions, one determines the primes contributing to the factorisation of $J(D_1, D_2)^2$. A delicate calculation using Deuring's theory of endomorphisms of elliptic curves in finite characteristic is necessary to determine the precise exponent in the factorisation formula. \\
	
	\par The first step of the analytic proof consists in producing a reformulation of \eqref{E: factorisation} more suitable to analytic manipulation. Let $F=\Q(\sqrt D)$ and consider the diagram of quadratic extensions		\begin{center}
		\begin{tikzpicture}[node distance = 1.2 cm, auto] \label{D: biquadratic Crelle}
		\node (Q) {$\mathbb{Q}$};
		\node (K1) [above of=Q, left of=Q] {$\Q(\sqrt{D_1})$};
		\node (K2) [above of=Q, right of=Q] {$\Q(\sqrt{D_2})$};
		\node (L) [above of=Q, node distance = 2.4 cm] {$\Q(\sqrt{D_1},\sqrt{D_2})$};
		\node (F) [above of=Q]{$F$};
		\draw[-] (Q) to node {} (K1);
		\draw[-] (Q) to node {} (K2); 
		\draw[-] (Q) to node {} (F);
		\draw[-] (K1) to node {} (L);
		\draw[-] (K2) to node {} (L);
		\draw[-] (F) to node {} (L);
		\end{tikzpicture}
	\end{center}
	where $\Q(\sqrt D_1, \sqrt D_2)$ is an unramified extension of the real quadratic field $F$. By Artin reciprocity, the genus character cutting out this extension can be viewed as a quadratic character of the narrow class group $\mathrm{Cl}^+(\cO_F)$ of $F$, which we denote by $\chi$. An elementary calculation shows that the explicit factorisation formula for $J(D_1,D_2)$ is equivalent to
	
	\begin{equation}\label{E: factorisation analytic}
	- \log \vert J(D_1,D_2) \vert^2= \sum_{\nu \in {\mathfrak d_{+}^{-1}}, \mathrm{Tr} \nu=1} \sum_{\mathfrak{n} \mid (\nu \mathfrak d)}\chi(n) \log \Nm(\mathfrak{n}),
	\end{equation}	
where $\mathfrak d^{-1}_+$ denotes the totally positive elements in the inverse different $\mathfrak d^{-1}$ of $\mathcal{O}_F$. The right hand side of \eqref{E: factorisation analytic} evokes formulae appearing in the work of Siegel on special values of Dedekind $\zeta$-functions for totally real fields \cite[\S 2]{Siegel}.  
	For $\nu \in F$, denote $(\nu_1, \nu_2)$ its real embeddings.  For $\mathbf z=(z_1,z_2) \in \cH^2_{\infty}$ and $(\nu, \mu)$ in $F^2$, we let 
	\[
	\mathrm{Nm}(\mu \mathbf z+\nu):=(\mu_1z_1+\nu_1)\cdot (\mu_2z_2+\nu_2).	\]
	The formulae appearing in \emph{loc.~cit.} arise as restrictions of Hilbert Eisenstein series of weight $(k,k)$, where $k\geq 2$ is an even integer, and level $\SL_2(\cO_F)$, which are given by
\begin{equation} \label{E: Eisenstein Siegel}
	E_{F,k}(\mathbf{z})=\sum_{\mathfrak{a} \in \mathrm{Cl}(\mathcal{O}_F)} \mathrm{Nm}(\mathfrak{a})^k \sum_{(m,n) \in (\af^2)'/ \mathcal{O}^{\times}_F} \mathrm {Nm}(m\mathbf z+n)^{-k}, \qquad \text{for } \mathbf z \in \cH^2_{\infty},
 \end{equation}
 along the diagonal embedding $\Delta \colon \cH_{\infty} \to \cH^2_{\infty}$.
	 The diagonal restriction inherits modular properties from those of $E_{F,k}$. It is an elliptic modular form of weight $2k$ and level 1.	
	 Siegel's formulae are obtained by expressing the resulting form, whose constant term in  the Fourier expansion encodes the Dedekind special value,  in terms of a basis of the space of elliptic modular forms for small values of $k$. 
	  
	The strategy of Gross and Zagier consists of mimicking Siegel's approach in the degenerate case in which $k=1$. The series in (\ref{E: Eisenstein Siegel}) would not converge for $k=1$. The problem can be obviated by considering the family of real analytic Hilbert Eisenstein series of parallel weight one 
\begin{equation} \label{E: Eisenstein GZ Crelle}	
	\mathcal {E}_{F,1, s}(\mathbf z)=\sum_{\af \in \mathrm{Cl}(\mathcal{O}_F)} 
	\chi(\af) \Nm(\af)^{1+2s}  
	\sum_{(m,n) \in (\af^2)'/ \mathcal{O}^{\times} }  
	\mathrm {Nm}(m\mathbf z+n)^{-k}
	\left |\mathrm{Nm}(m\mathbf z+n)\right |^{-s}
	\mathbf y^s
\end{equation}
for a complex variable $s$ with $\mathrm{Re}(s)>0$. For fixed values $\mathbf z \in \mathcal{H}_{\infty}^2$, this defines a holomorphic function of $s$, and it can be extended meromorphically to all of $\C$. The strategy for constructing holomorphic weight 1 Eisenstein series, also known as \emph{Hecke's trick}, consists in taking the limit of $\mathcal E_{F,1,s}$ as $s$ tends to 0. 
In this setting, Hecke's trick fails at producing an interesting weight one Hilbert Eisenstein series: an unexpected cancellation occurs when computing the corresponding Fourier expansion. Gross--Zagier are drawn to consider instead the derivative $\mathcal E'_{F,1,s}$ at $s=0$. This is a real analytic weight 1 Hilbert modular form; it can be  restricted along the diagonal embedding $\Delta$ to obtain a (non-holomorphic) elliptic modular form of weight 2, as in \cite{Siegel}.  
From this real analytic form $\Delta^* \mathcal E'_{F,1}$ one can extract a holomorphic one by applying a projector onto the holomorphic subspace as in \cite{Sturm}.  The resulting modular form, denoted by $(\Delta^* \mathcal E'_{F,1})_{\hol}$, is the object of the following theorem, summarising the calculations in \cite[\S 7]{GZ1985}. 

\begin{theorem}\label{T: a1 GZ Crelle} The first Fourier coefficient $a_1$ of the holomorphic modular form $(\Delta^* \mathcal E'_{F,1})_{\hol}$  of weight 2 and level $\SL_2(\Z)$ satisfies the equation  
\[
a_1 \cdot \lambda =  \log \vert J(D_1,D_2) \vert^2+ \sum_{\nu \in {\mathfrak d_{+}^{-1}}, \mathrm{Tr} \nu=1} \sum_{\mathfrak{n} \mid (\nu \mathfrak d)}\chi(n) \log \Nm(\mathfrak{n}),
\]
for some $\lambda \in \C^{\times}$. 
\end{theorem}
Now as $(\Delta^* \mathcal E'_{F,1})_{\hol}$ is of weight $2$ and level $\SL_2(\Z)$, one abstractly knows that it is the zero function, hence $a_1=0$ and the equality (\ref{E: factorisation analytic}) follows.

\subsection{Gross--Zagier's work on BSD in analytic rank 1} \label{SS: GZ Inventiones}
Beyond the factorisation of differences of singular moduli, the circle of ideas of \cite{GZ1985} yielded some striking developments towards the Birch and Swinnerton-Dyer Conjecture.

Let $E/\Q$ be an elliptic curve. The $L$-function $L(E, s)$ is defined as a convergent infinite product for $\mathrm{Re}(s)>3/2$. If $E$ is modular, that is $L(f,s)=L(E,s)$ for a suitable modular form $f$, the $L$-function $L(E,s)$ admits analytic continuation to all of $\C$ and functional equation with centre $s=1$. The Birch and Swinnerton-Dyer Conjecture states that the algebraic rank, that is the rank of the group $E(\Q)$, is equal to the analytic rank, the order of vanishing of $L(E,s)$ at $s=1$.  

The main result in \cite{GZ1986} is the following. 

\begin{theorem}\label{T: GZ Inventiones} \cite[Thm. 7.3]{GZ1986} Suppose that $L(E,1)=0$. There exists a point $P \in E(\Q)$ such that
\[
L'(E,1) \doteq \Omega_E \langle P,  P\rangle\] 
where $\Omega_E$ is the real period of a regular differential on $E$,  the pairing $\langle \cdot , \cdot \rangle$  denotes the N\'eron-Tate height pairing on $E(\Q) \otimes \R$ and the equality denoted by $\doteq$ is defined up to scalars in $\Q^{\times}$. 
\end{theorem}

In particular, if the analytic rank is one, the rank of the Mordell--Weil group $E(\Q)$ is positive. This result was later complemented by Kolyvagin's work bounding the algebraic rank in analytic rank 0 and 1.
The point $P$ arises as a Heegner construction (see below).
Combining the results of Gross--Zagier and Kolyvagin, one can conclude that this Heegner construction produces non-torsion points in $E(\Q)$ precisely when the algebraic (or analytic) rank is 1.

\par Let us discuss some ingredients of the proof. Let $N$ be the conductor of $E$. The assumption that $E$ is modular implies that up to isogeny $E$ is a direct factor  of the Jacobian $J_0(N)$ of the modular curve of level $X_0(N)$. Let $D<0$ be a fundamental discriminant coprime to $N$, and let $K=\Q(\sqrt D)$.
Let $\tau \in \mathcal{H}_{\infty}$ be an element with discriminant $D$ defining a CM point of the modular curve of level $\Gamma_0(N)$. Then the divisor $c_{\tau}:=(\tau)-(\infty)$ defines an $H$-rational point of $J_0(N)$, where $H$ is the Hilbert class field of $K$. By Shimura reciprocity, the divisor
\[
c_D=\sum_{[\tau] \colon \mathrm{disc}(\tau)=D} c_{\tau},\]
where the sum runs over a set of representatives of the $\Gamma_0(N)$-orbits of points of discriminant $D$ in $\cH_{\infty}$, is $K$-rational. The point $P$ appearing in the statement of Theorem \ref{T: GZ Inventiones} is constructed as a suitable trace of the divisor $c_D$.

The vector space $J_0(N)(H) \otimes \R$ is endowed with a canonical height pairing, which can be described as a sum of local heights and is equivariant for the Hecke action. 

The formal $q$-expansion 
\begin{equation}\label{E: modular generating series}
G_D(q)=\sum_{n \geq 1} \langle  c_D, T_nc_D \rangle q^n 
\end{equation}
is the $q$-expansion of a modular forms of weight 2 and level $\Gamma_0(N)$. This  formal series $G_D(q)$ is a prototypical example of a \emph{modular generating series}. It packages pairings of Hecke translates of Heegner objects into a formal $q$-expansion. Its modularity follows from the fact that the Hecke action on the Jacobian of the modular curve factors through the Hecke algebra of modular forms of weight 2 and level $\Gamma_0(N)$.

The proof of Theorem \ref {T: GZ Inventiones} hinges on providing an alternative construction of the $f$-isotypic component of the series  $G_D$ which can be related to $L$-values. For this, the key input is the Rankin--Selberg method. It allows to reinterpret the derivative $L'(f/K,s)$ 
at $s=1$ as the $f$-isotypic component of the holomorphic projection of the product of 
\begin{itemize}
\item  a weight 1 theta series $\theta_1$ attached to $K$ and 
\item the derivative at $s=0$ of a family of real analytic Eisenstein series $\mathcal{E}_{1,s}$ of constant weight 1, parametrised by a complex variable $s$. 
\end{itemize}
The comparison between the modular generating series $G_D$ and the holomorphic projection of  $\theta_1 \cdot \mathcal {E}'$ involves matching local contributions for both. In particular, it requires computing local height pairings for the divisor $c_D$ at non-archimedean places via intersection theory, and at archimedean places via Green functions.

\subsection{Biquadratic extensions and modular generating series}
\label{SS: biquadratic framework}
 While the arithmetic applications of \cite{GZ1985} and \cite{GZ1986} are on the surface fairly different, the strategies of proof fit into a general framework.  The diagram \ref{D: biquadratic Crelle} can be viewed as a special case of a diagram of quadratic field extensions (or, more generally, of \'etale algebras of degree 2) of the form
	\begin{center}
		\begin{tikzpicture}[node distance = 1cm, auto]
		\node (Q) {$\mathbb{Q}$};
		\node (K1) [above of=Q, left of=Q] {$K_1$};
		\node (K2) [above of=Q, right of=Q] {$K_2$};
		\node (L) [above of=Q, node distance = 2cm] {$L$};
		\node (F) [above of=Q]{$F$};
		\draw[-] (Q) to node {} (K1);
		\draw[-] (Q) to node {} (K2); 
		\draw[-] (Q) to node {} (F);
		\draw[-] (K1) to node {} (L);
		\draw[-] (K2) to node {} (L);
		\draw[-] (F) to node {} (L);
		\end{tikzpicture}
	\end{center}
where $F$ is a real quadratic field (or a split quadratic algebra over $\Q$, that is $F \simeq \Q \times \Q$) and $K_i=\Q(\sqrt{D_i})$ are imaginary, so that $L$ is forced to be a CM extension of $F$. 

For the imaginary quadratic fields $K_1$ and $K_2$, one can produce  Heegner divisors $c_{D_i}$ attached to points of discriminant $D_i$ in the Jacobian of modular curves of suitable level. As in \eqref{E: modular generating series}, the global height pairings on the Hecke translates of $c_{D_i}$ can be parlayed into a modular generating series $G_{D_1, D_2}(q)=\sum_{n \geq 1} \langle c_{D_1}, T_n c_{D_2} \rangle q^n $. The setting of Gross--Zagier's work on BSD in analytic rank 1 corresponds to the degenerate case in which the imaginary quadratic fields are equal, $F$ is the split quadratic algebra over $\Q$,  and $L =K_1 \times K_1$.

Obtaining an explicit characterisation of the Fourier coefficients of $G_{D_1, D_2}$  is an essential step towards  the desired arithmetic applications: 
\begin{itemize}
\item In \cite{GZ1985}, the quantity $J(D_1, D_2)$ can be read off the archimedean contribution of the height pairing $\langle c_{D_1}, c_ {D_2} \rangle$ for the modular curve of level 1;
\item In \cite{GZ1986}, for a modular elliptic curve $E$ corresponding to a newform $f$, the $f$-isotypic component of $G_{D_1, D_2}$ encodes the height of the Heegner point appearing in the statement of \ref {T: GZ Inventiones}. 
\end{itemize} 	

The crucial idea in both \cite{GZ1985} and \cite{GZ1986} is to express the modular generating series $G_{D_1, D_2}$ in an alternative form by exploiting the derivative of a suitable family $\cG_s$ of real analytic Hilbert modular forms of constant parallel weight 1 for the quadratic $\Q$-algebra $F$, parametrised by the variable $s \in \C$. 
In \cite{GZ1985}, the role of $\cG_s$ is the Eisenstein series $E_{F,1,s}$. In the degenerate setting of \cite{GZ1986}, the role of $\cG_s$ is played by the pair of elliptic modular forms $(\mathcal{E}_{1,s}, \theta_1)$ viewed as a pair of functions on two copies of $\cH_{\infty}$ invariant under the action of a congruence subgroup of $\GL_2(\Q \times \Q)$. 

The modular generating series $G_{D_1, D_2}$ is then shown to be equal to the formal $q$-expansion of $(\Delta^*\cG')_{\hol}$, the elliptic (holomorphic) modular form obtained by: 
\begin{enumerate}[label=\roman*)]
	\item computing the \emph{derivative} $\cG'$ of $\cG_s$ at $s=0$. This yields a real analytic Hilbert modular form of parallel weight 1;
	
	\item pulling back $\cG'$ under the diagonal embedding $\Delta \colon \cH_{\infty} \to \cH_{\infty}^2$. The resulting \emph{diagonal restriction} $\Delta^* \cG'$ is a real analytic elliptic weight 2 modular form; 

	\item applying a \emph{holomorphic projection}  to $\Delta^* \cG'$. 
\end{enumerate}
The resulting holomorphic projection $(\Delta^* \cG')_{\hol}$ is finally  shown to equal the modular generating series $G_{D_1, D_2}$.
	
	In \S \ref{S: p-adic families} the structure of these proofs will be mimicked in the $p$-adic setting, by replacing CM cycles with their real quadratic counterparts, and real analytic families with $p$-adic families of modular forms.

\begin{remark}
In \cite{GZ1985}, the modular generating series $G_{D_1, D_2}$ does not appear explicitly. Because the calculation is carried out for the modular curve of level 1, which has genus 0,  such a generating series is identically zero. The calculation in the analytic proof of \emph{loc.~cit.} amounts to determining the degree 1 coefficient of the modular generating series as a sum of local contributions, thereby leading to the desired factorisation formula. 
\end{remark}

\begin{remark}
Note that in \cite{GZ1986}, the theta series $\theta_1$ appearing as the second component of the family $\cG_s$ should be thought of as constant in the variable $s$. This reflects the fact that in the archimedean setting, only real analytic Eisenstein series admit variations in families. This imposes some restrictions on the Heegner  cycles for which the Gross--Zagier strategy can be carried out.  In general, one would expect the relevant Hilbert family $\cG_s$ to deform a parallel weight 1 Hilbert theta series attached to the CM extension $L/F$. However, these theta series are usually cuspidal, and as such they do not admit archimedean deformations. By contrast, in the $p$-adic setting, cusp forms can often be $p$-adically deformed. This will allow additional flexibility in the settings arising in \S \ref{S: p-adic families}.
\end{remark}

	\section{$p$-adic families of modular forms and RM theory }\label{S: p-adic families}
As we have seen in the last section real analytic families of Eisenstein series play a vital role in the works of Gross and Zagier on CM theory. In the $p$-adic setting there is a direct analogue, the $p$-adic family of Eisenstein series and below we will review how it is used in RM theory. Before we go into details we provide a bit of context for the much richer theory of \emph{ $p$-adic} families of modular forms.

\subsection{A brief overview of $p$-adic families of modular forms} \label{SS: p-adic families intro}
The theory of $p$-adic families of modular forms (or more generally of automorphic forms) has its origins in the 70s when Serre observed that the Hecke eigenvalues of Eisenstein series can be $p$-adically interpolated and hence that  Eisenstein series can be viewed as a $p$-adic family parametrised by the weight. This simple yet striking observation combined with a general interest in establishing congruences of modular forms motivated the search for a theory of $p$-adic variations of modular forms. 

A $p$-adic modular form (of level $\mathrm{SL}_2(\Z)$) as defined by Serre in \cite{Serre73} is a power series 
$f = \sum_{n \geq 0} a_n q^n \in  \Q_p[[q]]$
such that there is a sequence of classical modular forms $f_i=\sum_{n \geq 0} a_{i,n} q^n $ that converges to $f$ (uniformly in the coefficients). 
There is an alternative geometric definition due to Katz \cite{Katz}, which involves the 
(rigid analytic versions of) modular curves and generalises the classical definition of modular forms as sections of certain line bundles. 

 In a family of $p$-adic modular forms the weight is one of the crucial $p$-adically varying parameters. In fact this parameter can be viewed to vary in the so-called \emph{weight space}, which is a rigid analytic space $\mathcal{W}$ whose $\C_p$-points 
are given by  
\[
	\mathcal{W}(\C_p)=\Hom_{\mathrm{cts}} (\Z_p^{\times}, \C_p^\times)\ ,
	\]
the group of continuous $\C_p^\times$-valued characters of $\Z_p^{\times}$ and which identifies with a finite union of open unit discs. The weight $k \in \Z$ of a classical modular form is viewed as a $\Q_p$-point of $\mathcal{W}$ by identifying it with the character $z \mapsto z^{k-1}$, a general $p$-adic modular form then comes with a weight $\kappa \in \mathcal{W}(\C_p)$.

For building a good geometric theory of families the space of $p$-adic modular forms turns out to be too big, a problem that both Hida and Coleman managed to overcome using the Hecke operator $U_p$ which acts on the space of $p$-adic modular forms. In a nutshell, Hida studies the \emph{ordinary subspace}, i.e., the subspace of the space of $p$-adic modular forms on which $U_p$ acts invertibly. Coleman on the other hand works with certain subspaces of the space of $p$-adic forms, the so called  \emph{overconvergent} $p$-adic modular forms. Overconvergence was already studied by Katz and is an extra condition which decreases the size of these spaces (although the result is still an infinite-dimensional space). Now, on these spaces of overconvergent forms the operator $U_p$ acts as a compact operator and the complement of the kernel of $U_p$ is well behaved. Coleman managed to use this fact to build $p$-adic families of modular forms \cite{coleman}. 

From a geometric point of view the theory for $p$-adic families of modular forms then culminates in the construction of the so-called eigencurve by Coleman and Mazur (\cite{colemanmazur}). The eigencurve is a rigid analytic curve $\mathcal{C}$ which lives over $\mathcal{W}$ and whose points correspond to overconvergent $p$-adic modular Hecke eigenforms that are not in the kernel of the $U_p$-operator. The \emph{classical points}, i.e., the points corresponding to classical modular forms, form a Zariski-dense subset of $\mathcal{C}$. The global geometry of the Coleman--Mazur eigencurve remains a fascinating and challenging subject in current research (see \cite{LTXZ} for very recent progress on some of the original questions of Coleman and Mazur). 

Let us highlight one important feature of the theory of eigencurves, namely the connection to the theory of Galois representations. To a classical modular eigenform $f$ one can associate a 2-dimensional Galois representation $\rho_f$, which is characterised by identifying the Hecke eigenvalues with traces of Frobenius elements (\cite{Deligne}, \cite[Chapter 9]{DS}). As it turns out, one can interpolate these Galois representations (or rather, their traces and determinants) on the eigencurve to a family of pseudorepresentations. Moreover the infinitesimal neighbourhood of a classical point on the eigencurve has an interpretation in terms of a (suitably refined) Galois deformation space. This allows the usage of tools like Galois cohomology in the study of the local geometry of the eigencurve. 

Since the works of Serre, Katz, Coleman and Mazur, $p$-adic forms have been introduced for other reductive groups $G$ over a number field $F$ by different techniques. In particular the definition of overconvergent modular forms has been generalised to the setting of more general Shimura varieties. Furthermore we have constructions of so-called eigenvarieties generalising the Coleman--Mazur eigencurve. Let us emphasise (as this is relevant below) that we have a good theory of $p$-adic families of Hilbert modular forms at hand. 

\subsection{A modular generating series for RM values}
\label{SS: DPV1}
While a priori singular moduli for real quadratic fields do not enjoy the same algebraicity properties as their imaginary quadratic counterparts, they can be packaged into modular generating series, which one can hope to relate to derivatives of modular forms in the $p$-adic setting. In \cite{DPV1}, Darmon, Pozzi and Vonk imitate the strategy of \cite{GZ1985} outlined in \S \ref{SS: GZ Crelle} by replacing the real analytic family of Hilbert Eisenstein series appearing in \eqref{E: Eisenstein GZ Crelle} with a $p$-adic family parametrised by the weight, and relate it to a modular generating series for RM values of theta cocycles.  We will describe their result and compare  it with the formalism presented in \S \ref{SS: biquadratic framework}.\\

\par 
Let $F$ be a real quadratic field of discriminant $D$, let $\chi$ be an odd character of $\mathrm{Cl}^{+}(\mathcal{O}_F)$, and let $p$ be a prime unramified in $F$. Given an odd integer $k \in \Z_{\geq 1}$,    consider the parallel weight $k$ Hilbert Eisenstein series in the variable $\mathbf z=(z_1, z_2)$ in $\mathcal H^2_{\infty}$ with Fourier expansion
	\begin{equation} \label{E: Eisenstein DPV}
	E_{F, \chi, k}^{(p)}(\mathbf z)=L^{(p)}(F,\chi,1-k)+4\sum_{\nu \in \mathfrak d_{+}^{-1}} \sum_{\mathfrak n  \mid (\nu \mathfrak d), \ p \nmid \Nm(\mathfrak n )} \chi(\mathfrak n) \Nm (\mathfrak n)^{k-1} e^{2 \pi i \mathrm{Tr}(\nu \cdot \mathbf z)}
	\end{equation}
	where $\mathrm{Tr}(\nu \cdot \mathbf z):=\nu_1z_1+\nu_2z_2$ for the real embeddings $\nu_1,\nu_2$ of $\nu \in F$ and 
	\[
	L^{(p)}(F, \chi, s)=\sum_{\mathfrak n,\  p \nmid \mathrm {Nm}(\mathfrak n)} \chi(\mathfrak n) \mathrm {Nm}(\mathfrak n)^{s}. 
	\]
These modular forms are obtained as $p$-stabilisations of Hilbert modular forms of level $\SL_2(\cO_F)$ defined similarly to \eqref {E: Eisenstein Siegel} at primes above $p$ in $F$. The $p$-stabilisation procedure allows their Fourier coefficients to interpolate into continuous functions of the variable $k \in \mathcal {W}$, at the cost of introducing $p$ into the level. 
This can be verified directly for the Hecke coefficients attached to $\nu \in \mathfrak d^{-1}_+$, and can be deduced for the constant term  following Serre's approach to the analyticity of the Kubota--Leopoldt $p$-adic $\zeta$-function. We denote the family as $\mathcal E_{F, \chi, k}$, where the weight $k$ is thought of as a variable $k \in \mathcal W$, as explained in \S \ref{SS: p-adic families intro}. \\
\par 

Darmon, Pozzi and Vonk translate the steps outlined in \S \ref {SS: biquadratic framework} to the current setting. In analogy to Gross--Zagier, who consider  a family of constant weight 1, it is natural to analyse the specialisation of the family $\mathcal E_{F, \chi, k}$ at $k=1$, in view of applying a derivative.  While the Fourier expansion of \eqref{E: Eisenstein DPV} may not vanish at $k=1$, there is no harm in first considering the pullback of the family $\mathcal E_{F, \chi, k}$ along the diagonal embedding $\Delta \colon \cH_{\infty} \hookrightarrow \cH_{\infty}^2$. 
For every classical weight $k$, the diagonal restriction of the weight $k$-specialisation $E_{F, \chi, k}^{(p)}$ gives a classical elliptic modular form of weight $2k$, and its Fourier coefficients interpolate $p$-adically as functions of $k \in \cW$.
This yields a $p$-adic family of elliptic modular forms $\Delta^* \mathcal E_{F, \chi, k}$. Its behaviour at $k=1$ depends on the splitting of $p$ in $F$ (see \cite[Thm.\ A] {DPV1}).

\begin{itemize}
\item When $p$ splits in $F$, the specialisation of $\Delta^* \mathcal E_{F, \chi, k}$ need not vanish. The coefficients of its Fourier expansion can be expressed in terms of intersection pairings of the geodesic attached to real quadratic points of discriminant $D$ on the modular curve of level $\Gamma_0(p)$ and the path between $(0, \infty)$ on the modular curve of level $\Gamma_0(p)$. 

\item When $p$ is inert in $F$, the specialisation at $k=1$ vanishes identically. 
\end{itemize}
The split setting is rather classical in flavour. In the inert setting, the vanishing happens for similar reasons as in \cite{GZ1985}, as the weight 1 specialisation of $\Delta^* \mathcal E_{F, \chi, k}$  can itself be viewed as the $p$-stabilisation of the diagonal restriction of an Eisenstein series of weight 1 and trivial level, which must be identically 0. In the latter setting, it is tempting to consider the derivative $\Delta^* \mathcal E'_{F, \chi, k}$ of the above family at $k=1$. Exploiting the vanishing of $\Delta^* \mathcal E_{F, \chi, k}$ at $k=1$, one can show that this derivative is a $p$-adic (and, in fact, even overconvergent) modular form of weight 2 and trivial tame level \cite[Thm.\ 2.1]{DPV1}. Following the Gross--Zagier strategy, one wishes to tweak the resulting $p$-adic modular form and produce a classical one. This can be achieved by means of Hida's ordinary projector, constructed by taking limits of suitable iterates of the $U_p$ operator. We denote the resulting modular form by 
$(\Delta^* \mathcal E'_{F, \chi, k})_{\ord}$. It turns out to be a modular generating series for the RM values of the winding theta cocycle $J_{(0, \infty)}$ at real quadratic points of discriminant $D$.

\begin{theorem}\cite[Thm.\ B]{DPV1} \label{T: DPV1}
Suppose that $p$ is inert in $F$. The ordinary projection of the diagonal restriction of $\mathcal E'_{F, \chi, k}$ is a classical modular form of weight 2 and level $\Gamma_0(p)$ with $q$-expansion
\begin{equation*}
(\Delta^* \mathcal E'_{F, \chi, k})_{\ord}	
=L'_p(F, \chi,0) +\sum_{n \geq 1} q^n  
\sum_{[\tau],\  \disc \tau=D} \chi(\mathfrak c_{\tau}) \log_p \left( \Nm_{\Q_p}^{\Q_{p^2}} (T_n J_{(0, \infty)})[\tau] \right) 
\end{equation*}
where the sum runs over $\SL_2(\Z)$-classes of points of discriminant $D$,
$\mathfrak c_{\tau}$ denotes the fractional ideal $\Z+\tau \Z$ if $\tau-\tau'>0$ and $
\sqrt{D}(\Z+\tau \Z)$ otherwise, and $\Q_{p^2}$ is the unique unramified quadratic extension of $\Q_p$.
\end{theorem}

Note that the points $\tau$ in the above sum belong to $\cH_p$ precisely because $p$ is inert in $F$. Unlike in Thm.\ \ref{T: a1 GZ Crelle}, the constructed modular form has no a priori reason to vanish since it has level $\Gamma_0(p)$. 

\begin{remark}
The modular form in Theorem \ref{T: DPV1} is in general non-zero. Its constant coefficient is the leading term of a $p$-adic $L$-function which can be expressed as the $p$-adic logarithm of a certain global units, which will be the subject of Theorem \ref{T: DPV2}.

\end{remark}

\begin{remark}The values $J_{(0, \infty)}[\tau]$ of the winding theta cocycle at RM points $\tau$ in $\cH_p$ should be thought of as multiplicative analogues of the differences of singular moduli $j(\tau_1)-j(\tau_2)$ appearing in Gross--Zagier \cite{GZ1985}, where the pair $(0, \infty)$ corresponds to the split quadratic form $Q_{(0, \infty)}(X,Y)=XY$. However, this setting is degenerate, because the solutions to $Q_{(0, \infty)}$ lie at the boundary of the $p$-adic upper half plane. In particular, the corresponding theta cocycle is not expected to have algebraic values at real quadratic points in general. In the general framework described in \S \ref{SS: biquadratic framework}, this setting should correspond to the diagram of biquadratic $\Q$-algebras of the form
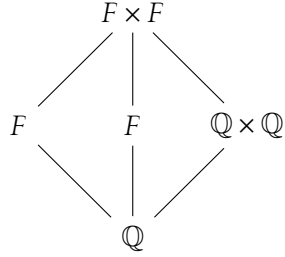
\begin{center}
		\begin{tikzpicture}[node distance = 1.5cm, auto]
		\node (Q) {$\mathbb{Q}$};
		\node (K1) [above of=Q, left of=Q] {$F$};
		\node (K2) [above of=Q, right of=Q] {$\Q \times \Q $};
		\node (L) [above of=Q, node distance = 3cm] {$F\times F$};
		\node (F) [above of=Q]{$F$};
		\draw[-] (Q) to node {} (K1);
		\draw[-] (Q) to node {} (K2); 
		\draw[-] (Q) to node {} (F);
		\draw[-] (K1) to node {} (L);
		\draw[-] (K2) to node {} (L);
		\draw[-] (F) to node {} (L);
		\end{tikzpicture}
	\end{center}
where the leftmost and rightmost algebras in the middle row  correspond to quadratic forms of discriminant $D$ and $Q_{(0, \infty)}$,  respectively.  The analogy with the settings of Gross--Zagier \cite{GZ1985} and \cite{GZ1986} is imperfect: in their work the archimedean place does not split in the quadratic fields $K_1$ and $K_2$.  
By contrast, in the setting of \cite{DPV1} the prime $p$ splits in the algebra $\Q \times \Q$. This is related to the fact that the Eisenstein series appearing  in \cite{DPV1} is attached to an arbitrary totally odd character of the class group $\mathrm{Cl}^+(\mathcal {O}_F)$, while in \cite{GZ1985} only genus characters are considered.
\end{remark}

\subsection{The values of the Dedekind--Rademacher theta cocycle} 
\label{SS: DPV2}

	In the archimedean setting, we discussed how modular generating series for Heegner cycles on the Jacobian of modular curves can be leveraged into arithmetic applications for the Heegner constructions discussed in \S \ref{S: Heegner CM}. Particularly, the main result in \cite{GZ1985} towards the Birch and Swinnerton-Dyer Conjecture is obtained by projecting a suitable modular generating series into the $f$-isotypic component for a newform $f$ corresponding to an elliptic curve via the Eichler--Shimura construction. \\
	Similarly, we will now discuss an arithmetic application of the ideas appearing in \S \ref{SS: DPV1} to the algebraicity of the RM values of the Dedekind--Rademacher theta cocycle $J_{\mathrm{DR}}$ defined in \S \ref{SS: Heegner objects RM}. The values  of the Dedekind--Rademacher cocycle at real multiplication points should be thought of as real quadratic analogues of the elliptic units, in light of Conjecture \ref {C: Darmon-Dasgupta}. 
	The cocycle $J_{\mathrm{DR}}$ is the Eisenstein class in the module of analytic theta cocycles $H^1(\Gamma, \mathcal A^{\times}/\C_p^{\times})$, as its definition involves the periods of the weight $2$ Eisenstein series of level $\Gamma_0(p)$. Accordingly, the main result of \cite{DPV2} is obtained by projecting a suitable modular generating series similar to the one appearing in Theorem \ref{T: DPV1} onto  its Eisenstein subspace. 
	
	However, a significant new element occurs in this setting that does not have an archimedean counterpart. The modular generating series in question is obtained by considering the derivative of a $p$-adic family of \emph{cusp forms} deforming over an appropriate weight space for Hilbert modular forms for a real quadratic field $F$. The coefficients of its derivative are richer than those of Eisenstein series and encode information about the arithmetic of the Hilbert class field of $F$. 
	
The main result \cite{DPV2} can be formulated as follows.  

\begin{theorem}\label{T: DPV2}
Let $D > 0$ be a fundamental discriminant, and  let $p \nmid D$ be a prime number. Let $\tau \in \cH_p$ be a point of discriminant $D$.  
Let $F = \Q(\tau)$ and let $H$ be the narrow Hilbert class field of $F$. 
The values of the Dedekind--Rademacher cocycle at $\tau$ satisfy

\begin{equation*}\label{E: values Dedekind-Rademacher}
J_\text{DR}[\tau] = u_\tau^{12} \text{ modulo }(p^\Z \times \text{ roots of unity in }\Q_{p^2}^\times)
\end{equation*}
for an element $u_\tau \in \oo_H[1/p]^\times\otimes\Q$.
\end{theorem}

Prior to this work, the equality of local norms
\begin{equation} \label {E: norms}
\Nm^{\Q_p^2}_{\Q_p} (J_\text{DR}(\tau) )\doteq \Nm^{\Q_p^2}_{\Q_p}(u_\tau^{12}) 
\end{equation}
up to $p$-powers and roots of unity in $\Q_{p^2}$ was already known, as both sides of  (\ref{E: norms}) had been related to the derivative of a $p$-adic $L$-function (for a definition, see \cite [\S 4]{DD06}). More precisely, Thm. 4.1 in \emph{loc.~cit.} relates the derivative of this $p$-adic $L$-function to the norm of the RM values of the Dedekind--Rademacher cocycle.
On the other hand, the $p$-unit $u_{\tau}$ is an example of the  Gross--Stark units, appearing in the conjectures of Stark and Gross about derivatives of Artin $L$-functions  and their $p$-adic analogues. In particular, the formula for the derivative of the $p$-adic zeta function above and the norm of the $p$-unit $u_{\tau}$ was proved by Darmon, Dasgupta and Pollack in their work on the Gross--Stark Conjecture. \\
The main contribution of the refinement of \eqref{E: norms} in Theorem \ref{T: DPV2} is proving the \emph{algebraicity} of the RM values of the Dedekind--Rademacher cocycle itself, for which the equality up to local norms would not suffice. This  establishes cases of Conjecture \ref{C: Darmon-Dasgupta}, providing theoretical evidence in favour of pursuing a theory of real multiplication via  $p$-adic methods.	\\

\begin{remark}
The unit $u_{\tau}$ in the statement of Thm.~\ref{T: DPV2} can be characterised more precisely. It is trivial if $F$ has a unit of negative norm, or equivalently, if the narrow class field of $H$ is totally real. Otherwise, it is a non-trivial element in the subspace of $\cO_{H}[1/p]^{\times} \otimes \Q$ on which the complex conjugation in $H$ acts as $-1$. The recent breakthroughs of Dasgupta and Kakde on the Brumer--Stark Conjecture imply that this unit in fact lies in  
$\cO_{H}[1/p]^{\times} \otimes \Z[1/2]$.
\end{remark}

\par 
The crucial intuition behind the strategy in \cite{DPV2} is that removing the norm in \eqref{E: norms} could be achieved by removing the norms in the modular generating series in Theorem \ref{T: DPV1} obtained from the derivative of the $p$-adic family  of Eisenstein series. This approach requires finding a suitable family of $p$-adic modular forms and exploiting its derivatives to produce the desired modular generating series. \\
Fortunately, the setting allows a lot of flexibility due to a lucky coincidence. The weight 1 Eisenstein series $E^{(p)} _{F, \chi, 1}$ considered in \eqref{E: Eisenstein DPV} happens to be cuspidal when viewed as an  overconvergent $p$-adic Hilbert modular form. In fact, it is an example of a critical $p$-stabilisation of a classical Eisenstein series (considered, for example by Bellaiche and Chenevier for elliptic modular forms \cite{BC}). As such, it defines a point on the eigenvariety parametrising $p$-adic families of Hilbert cusp forms for $F$. Unlike Eisenstein series, which only vary in parallel weight, Hilbert families of cusp forms range over a larger weight space parametrising pairs of weights $(k_1, k_2)$.\\
 Cuspidal variations of Hilbert modular forms are, in general, entirely inexplicit. However, their attached Galois representations provide a useful tool  to tackle a local description. In particular, they can be exploited to entirely determine the Fourier expansion of the derivative of a family $\mathcal F_{k, 2-k}$ specialising to $E^{(p)} _{F, \chi, 1}$ at $k=1$,
in the \emph{antiparallel} direction of weights $(k_1, k_2)$ satisfying 
$
k_1+k_2=2.
$
The diagonal restriction of the derivative $\mathcal F_{k, 2-k}$ with respect to $k$, denoted by $\Delta^*\mathcal F',$  is again an overconvergent $p$-adic modular form. Its ordinary projection $(\Delta^*\mathcal F')_{\ord}$, (up to an explicit correction term obtained from the diagonal restriction of explicit Eisenstein series) yields the following theorem, which is the key stepping towards Thm.\ \ref{T: DPV2}.
\newpage

\begin{theorem}\cite[Thm.~C]{DPV2} There is a classical modular form of weight 2 and level $\Gamma_0(p)$ with $q$-expansion
\[
\log_p(u_{\tau})+\sum_{n\geq 1}\log_p(T_nJ_{(0,\infty)}[\tau])q^n,
\]
for a unit $u_{\tau} \in \mathcal O_H[1/p]^{\times}\otimes \Q$. 
\end{theorem} 
From this result, the relation between the unit $u_{\tau}$ and the RM values of the Dedekind-Rademacher cocycle can easily be deduced by decomposing the winding cocycle $J_{(0, \infty)}$ with respect to an eigenbasis of the space of analytic theta cocycles, and in particular determining its projection onto the Eisenstein eigenspace spanned by $J_\mathrm{DR}$.

\begin{remark}
The relevance of cuspidal deformations of the Eisenstein series $E_{F, \chi, 1}^{(p)}$ towards Thm. \ref{T: DPV2} should not be surprising. A cuspidal variation of this Eisenstein series in the parallel weight direction  played a crucial role in the proof of the Gross--Stark Conjecture in \cite{DDP} implying the equality \eqref{E: norms}. This fits into a general theme relating congruences between cusp forms and Eisenstein series to $L$-functions, initiated by Ribet in \cite{Ribet} and often referred to as \emph{Ribet's method}. Theorem \ref{T: DPV2} is a special case of a far more general result of Dasgupta and Kakde  \cite{DK}, providing explicit $p$-adic formulae for units in abelian extensions of totally real fields. Their result can be seen as the  culmination of a program of Dasgupta and collaborators settling Hilbert's Twelfth problem for totally real fields via $p$-adic methods. Their approach makes use of generalisations of Ribet's method. However, unlike in \cite{DPV2}, the strategy is applied to variations of Hilbert modular forms over finite group rings, rather than over a $p$-adic weight space.
\end{remark}

	\bibliography{overviewRMC}
	\bibliographystyle{alpha}
	
\end{document}